\documentclass[11pt,twoside]{amsart}
\usepackage[dvipsnames]{xcolor}
\usepackage{hyperref}
\hypersetup{colorlinks=true, linkcolor=Blue, citecolor=Maroon}
\usepackage{amsthm, amsmath, amscd, amssymb,centernot,txfonts}
\usepackage{tikz-cd}
\usepackage{lipsum}
\usepackage{cancel} 
\usepackage{multicol}
\usepackage{amsfonts}
\usepackage{mathrsfs}
\usepackage[all]{xy}
\usepackage[left=2.5cm,top=2.5cm,bottom=3cm,right=2.5cm]{geometry}
\usepackage{dirtytalk}
\usepackage{mathtools}
\usetikzlibrary{shapes,arrows}
\usepackage{etoolbox}
\usepackage{dynkin-diagrams}
\usepackage{enumitem}
\usepackage{chngpage}
\usepackage{adjustbox}
\usepackage[normalem]{ulem}
\usepackage[utf8]{inputenc}
\usepackage{fourier} 
\usepackage{array}
\usepackage{makecell}
\usepackage{comment}
\usepackage{cancel}

\newcommand{\longsquiggly}{\xymatrix{{}\ar@{~>}[r]&{}}}

\theoremstyle{plain}
\numberwithin{equation}{section}
\newtheorem{theorem}{Theorem}[section]
\newtheorem{proposition}[theorem]{Proposition}
\newtheorem{lemma}[theorem]{Lemma}
\newtheorem{corollary}[theorem]{Corollary}

\theoremstyle{definition}
\newtheorem{remark}[theorem]{Remark}

\newtheorem{example}[theorem]{Example}

\newtheorem{set-up}[theorem]{Set-up}
\newtheorem{definition}[theorem]{Definition}

\newcommand*{\QEDA}{\hfill\ensuremath{\blacksquare}}

\tikzstyle{decision} = [diamond, draw, , 
    text width=4.5em, text badly centered, node distance=3cm, inner sep=0pt]
\tikzstyle{block} = [rectangle, draw, , 
    text width=10em, text centered, rounded corners, minimum height=2em]
\tikzstyle{block1} = [rectangle, draw, , 
    text width=5em, text centered, rounded corners, minimum height=2em]    
\tikzstyle{line} = [draw, -latex']
\tikzstyle{cloud} = [draw, ellipse,, node distance=3cm,
    minimum height=2em]
\begin{document}

\title[Projective smoothing of varieties with simple normal crossings]{Projective smoothing of varieties with simple normal crossings}

\author[P. Bangere]{Purnaprajna Bangere}
\address{Department of Mathematics, University of Kansas, Lawrence, USA}
\email{purna@ku.edu}

\author[F.J. Gallego]{Francisco Javier Gallego}
\address{Departamento de \'Algebra, Geometr\'ia y Topolog\'ia and Instituto de Matem\'atica Interdisciplinar,
Universidad Complutense de Madrid, Spain}
\email{gallego@mat.ucm.es}

\author[J. Mukherjee]{Jayan Mukherjee}
\address{Department of Mathematics, Oklahoma State University, Stillwater, USA}
\email{jayan.mukherjee@okstate.edu}


\maketitle
\begin{abstract}
In this article, we introduce a new approach to show the existence and smoothing of simple normal crossing varieties in a given projective space. Our approach relates the above to the existence of nowhere reduced schemes called ribbons and their smoothings via deformation theory of morphisms. As a consequence, we prove results on the existence and smoothing of snc subvarieties $V \subset \mathbb{P}^N$, with two irreducible components, each of which are Fano varieties of dimension $n>2$, embedded inside $\mathbb{P}^{N}$ for effective values of $N$, by the complete linear series of a line bundle $H$. The general fibers of the resulting one parameter families are either smooth Fano, Calabi-Yau or varieties of general type, depending on the positivity of the canonical divisor of their intersections. An interesting consequence of projective smoothing is that it automatically gives a smoothing of the semi-log-canonical (slc) pair $(V, \Delta)$, where $\Delta = cH$, $c < 1$, is a rational multiple of a general hyperplane section of $H$. For threefolds, we are able to give explicit descriptions of the smoothable snc subvarieties due to the classification results of Iskovskikh-Mori-Mukai. In particular, we show the existence of unions V = $Y_1 \bigcup_D Y_2 \subset \mathbb{P}^N$, where $Y_i$'s are smooth anticanonically (resp. bi-anticanonically) embedded Fano threefolds, intersecting along $D$, where $D$ is either a del-Pezzo surface or a $K3$ surface (resp. a smooth surface with ample canonical bundle) and their smoothing in $\mathbb{P}^N$ to smooth Fano or Calabi-Yau threefolds (resp. to threefolds with ample canonical bundle) for various values of $N$ between $10$ and $163$. In cases when the general fiber is a smooth Fano or Calabi-Yau threefold, one can choose $c$ such that $(V, \Delta)$ is a Calabi-Yau pair while in all cases $c$ can be chosen so that $(V, \Delta)$ is a stable pair.

\end{abstract}

\section{Introduction}

The topic of degenerations of algebraic varieties to varieties with simple normal crossings has been a topic of much interest in algebraic geometry. The topic of abstract smoothing of normal crossing varieties, particularly the ones having invariants that of either Fano or Calabi-Yau have been studied in \cite{Ku77}, \cite{F83}, \cite{KN94}, \cite{TN10}, \cite{TN15}, \cite{FFR21}. In this article we study projective degenerations of algebraic varieties to varieties with simple normal crossings. In addition to Fano and Calabi-Yau type, a large class of examples of smoothable simple normal crossings in this article are of general type. There are multiple motivations and implications for pursuing the topic of projective degenerations apart from the fact that it is an interesting question in its own right. For example, some of the motivations include compactifications of moduli of polarized varieties since snc varieties have semi-log-canonical singularities  (see \cite{Vie95}, \cite{KSBA88}, \cite{KX20}), the study of Gaussian maps (\cite{CLM93}), semistability of $q-$th Hilbert points for low values of $q$ (see \cite{Fed18}, \cite{DFS16}, \cite{HH09}, \cite{HH13}), syzygies of algebraic varieties (see \cite{BE95}, \cite{Deo18}, \cite{RS22}) among others. 

\smallskip
\textcolor{black}{An important class of varieties in the context of higher dimensional geometry is the class of Fano varieties. A lot is known about their geometry, boundedness etc in any dimension (see \cite{KMM}, \cite{Bir21}) including a complete classification in dimension three in the work of Iskovskikh-Mori-Mukai  \cite{I77}, \cite{I78}, \cite{MM82} and they are reasonably well-understood. Hence a motivating topic to pursue is the study of degenerations of smooth varieties into snc varieties whose irreducible components are Fano varieties. So we start by asking the following natural questions arising in this context, answers to which have important consequences as clarified below.}

\begin{itemize}
    \item[(1)] Fix a positive integer $N$. Does there exist a simple normal crossing variety $V$ which is a union $Y_1 \bigcup_D Y_2 \subset \mathbb{P}^N$ of smooth $n-$ anticanonically embedded Fano varieties $Y_i$ that intersect along a smooth divisor $D$ which is either
    \begin{itemize}
        \item[(a)] a Fano variety 
        \item[(b)] a Calabi-Yau variety or
        \item[(c)] a variety with ample canonical bundle ?
       \end{itemize}
\end{itemize}

Such a union gives a Gorenstein simple normal crossing variety in $\mathbb{P}^N$ for effective values of $N$ with invariants that of a Fano variety, a Calabi-Yau variety and a variety with ample canonical bundle in cases (a), (b) and (c) respectively. Further, pulling back a suitable rational multiple $\Delta = cH$ with $c \leq 1$ of a general hyperplane section $H$, the pair $(Y_1 \bigcup_D Y_2, \Delta)$ gives rise to slc Calabi-Yau pairs in cases (a) and (b) and stable pairs in cases (a), (b), (c). It is to be noted that the existence of such a union is not just about constructing a very ample line bundle on $Y_1 \bigcup_D Y_2$ but is also about showing that after a projection to a suitable linear subspace it remains an embedding. This is due to the fact that a very ample line bundle $L$ on an abstract union, embeds each component $Y_i$ degenerately inside $\mathbb{P}(H^0(L))$ in general.

\smallskip
 Once the existence of such an embedded simple normal crossing variety is established, the next important question is:
\begin{itemize}
    \item[(2)]  \textcolor{black}{Are these} simple normal crossing varieties $V = Y_1 \bigcup_D Y_2$ smoothable inside $\mathbb{P}^N$ ? 
\end{itemize}
    
It is important to note here that an abstract smoothing of the simple normal crossing in general will not be achieved inside the given projective space. The relevance of projective smoothing (apart from its many other applications as we mentioned above) is in the fact that it provides a smoothing of not only $Y_1 \bigcup_D Y_2$ but also of the \textcolor{black}{log pair} $(Y_1 \bigcup_D Y_2, \Delta)$. \par
We introduce an \textcolor{black}{algorithmic process} for a \textcolor{black}{systematic} construction and smoothing of simple normal crossing varieties consisting of two irreducible components meeting along a smooth divisor inside a projective space and in particular answer the above questions in the next theorem and the corollaries that follow:

\smallskip

\begin{theorem}\label{Fano introduction}(see Theorem \ref{Fano}) Let $Y$ be a smooth, projective Fano variety of dimension $d \geq 3$, let $H$ be a very ample line bundle on $Y$, let  $N=h^0(H)-1$ and, abusing the notation, call $Y$ the image in $\mathbb P^N$ of the embedding induced by the complete linear series $|H|$. Let $L$ be a line bundle such that $L^{-1}$ is ample and there exist a smooth divisor  $D$ in  $|L^{-1}|$ and a smooth divisor  in $|L^{-2}|$. Assume the following conditions hold: 

\begin{enumerate}
    
     \item[(a)] $H \otimes L$ is base point free.
     \item[(b)] $N \geq 2d+1$. 
     \item[(c)] One of the following conditions hold:
     \begin{enumerate}
        \item[(i)] $\omega_Y^{-1} \otimes L$ is ample; or
        \item[(ii)] $L =  \omega_Y$ and $h^{1,1}(Y) = 1$ when $d = 3$ or
        $L= \omega_Y$ and $h^{1,d-2}(Y) = 0$
       when $d \geq 4$; 
        \item[(iii)]  $L^{-1} \otimes  \omega_Y $ is ample,  
        $H^1(H  \otimes L)=0$
        and $Y$ satisfies Bott vanishing theorem. 
    \end{enumerate}
\end{enumerate}

Then
\begin{itemize}
    \item[(1)] There exist simple normal crossing varieties $V = Y_1 \bigcup_D Y_2$, embedded inside $\mathbb{P}^N$ with two irreducible components $Y_i$, such that $Y_1$ is the subvariety $Y$ and $Y_2$ is a deformation of the subvariety $Y$ inside $\mathbb{P}^N$ intersecting along the subvariety $D$.
    \item[(2)] The snc subvarieties $V$ are smoothable inside $\mathbb{P}^N$ into smooth Fano varieties or Calabi-Yau varieties or varieties of general type in cases $(c)$(i), $(c)$(ii) and $(c)$(iii) respectively. 
    \item[(3)] The locally trivial deformations of $V$ inside $\mathbb{P}^N$ are unobstructed
    and the first order locally trivial deformations form a subspace in the tangent space of the Hilbert scheme at $[V]$ of codimension $h^0(L^{-2}|_D)$. 
    Furthermore, this subspace is the tangent space at $[V]$ of an irreducible locus of the Hilbert scheme, which has codimension $h^0(L^{-2}|_D)$ and is smooth at $[V]$. The singularities of the subschemes parameterized by this locus are normal crossing singularities;  in particular, they are non-normal and semi-log-canonical. They are analytically isomorphic to $(x_1^2+x_2^2 = 0 \subset \mathbb{C}^N)$. 
    \item[(4)] The general fiber of any one-parameter family of deformations of $V$ whose image under the Kodaira-Spencer map is a first order locally trivial deformation of $V$, is once again a simple normal crossing variety $V' = Y_1' \bigcup_D' Y_2'$, where both $Y_1'$ and $Y_2'$ are deformations of the subvariety $Y$ and $D'$ is a deformation of the subvariety $D$. 

\end{itemize}

\end{theorem}

To understand the algorithm mentioned before the Theorem above, we state the definition of nowhere reduced scheme structures called ribbons.
\begin{definition}\label{def of a ribbon}(see \cite{BE95})
    A ribbon on a reduced connected scheme $Y$ is a nowhere reduced scheme $\widetilde{Y}$ with a closed embedding $Y \hookrightarrow \widetilde{Y}$ such that the ideal sheaf $\mathcal{I}$ of $Y$ inside $\widetilde{Y}$ satisfies (i) $\mathcal{I}^2 = (0)$ (ii) $\mathcal{I}/\mathcal{I}^2$ as a module over $\mathcal{O}_Y$ is locally free of rank one. The line bundle $L = \mathcal{I}/\mathcal{I}^2$ is called the conormal bundle of the ribbon.
\end{definition}
The algorithm presented here is based upon the following strategy (see Theorem \ref{main}, Theorem \ref{intrinsic})  :
\begin{itemize}
    \item[(a)] Start with an embedding of a smooth variety $Y \xhookrightarrow{i} \mathbb{P}^N$ and a smooth divisor $D$ in $Y$.
    \item[(b)] Show the existence of an embedded ribbon $\widetilde{Y} \hookrightarrow \mathbb{P}^N$ on $Y$ with conormal bundle $\mathcal{O}_Y(-D)$.
    \item[(c)] Construct a flat family of simple normal crossings which degenerates to the ribbon $\widetilde{Y}$.
    \item[(d)] Show that $\widetilde{Y}$ has an embedded smoothing   
    \item[(e)] $\widetilde{Y}$ represents a smooth point of the corresponding Hilbert scheme.
\end{itemize}

Step (c) is one of the key steps in the process. It is shown in this article that blowing up the embedded ribbon $\widetilde{Y}$ along the reduced Weil divisor $D$ gives us a first order deformation of $Y \hookrightarrow \mathbb{P}^N$ fixing $D$ (see Remark \ref{blow-up of a ribbon}). Extending that deformation to a one-parameter family $T$, and taking the union with the trivial embedded deformation of $Y \hookrightarrow \mathbb{P}^N$ over $T$, is the sketch of proof of Step (c). Some of the technical issues are to show that the central fibre of the family we get is indeed the ribbon $\widetilde{Y}$ and that the general fibres are snc. Step (d) ensures that the snc's are further smoothable. The methods employed in Step (d) give us a good understanding of the general smooth fibre, for example their Picard groups and the relation between the generator of the Picard group and the hyperplane section of the embedding. The smoothable simple normal crossings that we construct are general in the locus of their locally-trivial deformations.

Specializing to threefolds, we use the classification of smooth Fano threefolds by Iskovskikh-Mori-Mukai (see \cite{I77}, \cite{I78}, \cite{MM82}) to give a list of unions of anticanonically embedded Fano-threefolds that intersect along del-Pezzo surfaces and $K3$ surfaces and a list of unions of bi-anticanonically embedded Fano-threefolds that intersect along surfaces with ample canonical bundle all of which are projectively smoothable to smooth Fano, Calabi-Yau and threefolds with ample canonical bundle respectively. As mentioned earlier, this gives a larger list of smoothable polarized Calabi-Yau pairs and stable pairs. We denote the various deformation types of Fano threefolds using the notation in the website \cite{Bel24}, \url{https://www.fanography.info/}.  \par

\medskip

\begin{corollary}(see Corollary \ref{Fano-threefolds}, Corollary \ref{log CY and gen type pairs from Fano}, Example \ref{examples of log pairs from Fano}): {\bf{\underline{Smoothable union of anticanonically}}} \\{\bf{\underline{embedded Fano-threefolds intersecting along del-Pezzo surfaces}}} Applying Theorem \ref{Fano introduction}, $(c)$(i), we show that for each deformation type of Fano threefolds with index at least two and $|\omega_Y^{-1}|$ very ample, i.e, one of families $1.12(b)-1.16$ and $2.32,2.35,3.27$, there exists one parameter families $\mathcal{X} \hookrightarrow \mathbb{P}^{N}_T$ over $T$, such that $\mathcal{X}_0 = Y_1 \bigcup_D Y_2$ is a snc variety where $Y_i \hookrightarrow \mathbb{P}^{N}$, $i=1,2$ are anticanonically embedded Fano threefolds of the same deformation type, and $D = Y_1 \bigcap Y_2$ is a smooth sub anti-canonical del-Pezzo surface, while the general fibre $\mathcal{X}_t$ is a smooth Fano threefold embedded by a proper sublinear series of $|-2K_{\mathcal{X}_t}|$ in all cases other than $1.15$, $|-3K_{\mathcal{X}_t}|$ for family $1.15(a)$, $|-\displaystyle\frac{3}{2}K_{\mathcal{X}_t}|$ for family $1.15(b)$. In each case we identify what the smooth fibres are once again according to the classification. In each case, we list the possible choices of $c$ such that $(Y_1 \bigcup_D Y_2, \Delta = cH)$, is a Calabi-Yau pair or a stable pair (which are hence smoothable). 

\end{corollary}
\medskip

\begin{corollary}(see Corollary \ref{Calabi-Yau threefolds}): {\bf{\underline{Smoothable union of anticanonically embedded Fano-threefolds}}} \\ {\bf{\underline{intersecting along $K3$ surfaces}}} Applying Theorem \ref{Fano introduction}, $(c)$(ii), we show that for each deformation type $1.5-1.10$ and $1.12-1.17$ (resp. $1.1-1.4$ and $1.11$), each of which has Picard rank one, there exists one parameter families $\mathcal{X} \hookrightarrow \mathbb{P}^{N}_T$ over $T$, such that $\mathcal{X}_0 = Y_1 \bigcup_D Y_2$ is a snc variety where $Y_i \hookrightarrow \mathbb{P}^{N}$, $i=1,2$ are anticanonically (resp. bi-anticanonically) embedded Fano threefolds of the same deformation type, and $D = Y_1 \bigcap Y_2$ is a smooth anticanonical $K3$ surface,  while the general fibre $\mathcal{X}_t$ is a smooth Calabi-Yau threefold with Picard number $1$ which is not always prime. For families $1.5-1.10$, the fibre over the point $(D \hookrightarrow \mathbb{P}^N)$ of the flag Hilbert scheme parameterizing $(D \hookrightarrow Y \hookrightarrow \mathbb{P}^N)$ with varying $Y$ with the same Hilbert polynomial as $Y_i$'s is irreducible and consequently, any subvariety $V$ of the form $V=Y_1 \bigcup_D Y_2$ is smoothable. 
\end{corollary}
\medskip
\begin{corollary}(see Corollary \ref{Fano with Bott}):
{\bf{\underline{Smoothable union of bi-anticanonically embedded Fano-threefolds}}} \\ {\bf{\underline{intersecting along surfaces with ample canonical bundle}}} Applying Theorem \ref{Fano introduction}, $(c)$(iii), we show that for each deformation type of Fano threefolds that satisfy Bott vanishing theorem, which by \cite{Tot23} are families $1.17$, $2.33-2.36$, $3.25-3.31$, $4.9-4.12$, $5.2-5.3$ (toric), $2.26$, $2.30$, $3.15-3.16$, $3.18-3.24$, $4.3-4.8$, $5.1$, $6.1$ (non-toric), there exists one parameter families $\mathcal{X} \hookrightarrow \mathbb{P}^{N}_T$ over $T$, such that $\mathcal{X}_0 = Y_1 \bigcup_D Y_2$ is a snc variety where $Y_i \hookrightarrow \mathbb{P}^{N}$, $i=1,2$ are bi-anticanonically embedded Fano threefolds of the same deformation type, and $D = Y_1 \bigcap Y_2$ is a smooth surface with ample canonical bundle, while the general fibre $\mathcal{X}_t$ is a smooth threefold embedded by a sublinear series of the bicanonical system. The Picard number of $\mathcal{X}_t$ is equal to the Picard number of $Y_i$. 
\end{corollary}
\smallskip

\par In \cite{CLM93}, the authors show (\textcolor{black}{among other things}) degeneration of $K3$ surfaces of degree $2g-2$ into a union of two rational normal scrolls each of degree $g-1$ inside $\mathbb{P}^g$ meeting along an elliptic curve, anticanonical in each scroll. Our results on projective degenerations of Calabi-Yau varieties into a union of anticanonically embedded Fano varieties meeting along a $K3$ surface are higher dimensional analogues of the result in \cite{CLM93}. \textcolor{black}{But there is an interesting difference in the higher dimensional case. By the result of Kulikov (see \cite{Ku77}), in a type II degeneration of $K3$ surfaces with two irreducible components, both components must be rational, while in Theorem \ref{Calabi-Yau threefolds}, we see  examples of degenerations of polarized Calabi-Yau varieties into snc union of anticanonically embedded Fano varieties, neither of whose components are rational (see Remark \ref{deviation from Kulikov}).} In \cite{BGM24a}, \cite{BGM24b}, \cite{BGM24c}, we list some further work in progress.
\smallskip

We want to draw some contrasts and parallels with earlier results. Most of the previous work on smoothability of normal crossing varieties are about abstract smoothing of Fano or Calabi-Yau normal crossings. More specifically, the smoothing of $d-$ semistable Fano or Calabi-Yau normal crossings have been studied in great detail in \cite{Ku77}, \cite{F83}, \cite{KN94}, \cite{TN15}, while in \cite{TN10}, \cite{TN15} and \cite{FFR21}, more general results that do not require $d-$ semistability are proven. In contrast, the algorithmic process discussed above constructs projective smoothing of non $d-$ semistable simple normal crossings of different Kodaira dimension, for example, Fano, Calabi-Yau or general type by starting from the same embedding $Y \hookrightarrow \mathbb{P}^N$ and by varying the positivity of $D$. 
We also refer to \cite{Dr\'e15} for a related construction involving ribbons in dimension one, \cite{Fel23} for results on deformations of line bundles, \cite{CLM19} for different approaches to smoothing normal crossing varieties, \cite{CR23} for general construction of log structures on degenerate varieties and \cite{FFP23} for computing $T_V^1$ on semi-smooth varieties $V$.

\smallskip

Now we describe the structure of the paper. In Section \ref{main results} we develop the methods that are foundations for the rest of the paper. The main results of this section are Theorem \ref{main}, Theorem \ref{intrinsic} and Theorem \ref{theorem.irreducible}. In Theorem \ref{main}, we show how the existence of a ribbon $\widetilde{Y}$ on $Y$ embedded inside $\mathbb{P}^N$ gives rise to the existence of a one-parameter family of simple normal crossings that degenerate to the ribbon. Theorem \ref{intrinsic} gives a sufficient condition for an embedded smoothing of the simple normal crossings constructed in Theorem \ref{main}. Additionally in Theorem \ref{intrinsic}, we study the embedded locally-trivial deformations of the SNC subvarieties to understand the dimension of the locus inside the Hilbert scheme parameterizing subvarieties with the same singularity type. Theorem \ref{theorem.irreducible} connects all the above with the flag Hilbert schemes parameterizing $(D \hookrightarrow Y \hookrightarrow \mathbb{P}^N)$ where $D$ and $Y$ varies with fixed Hilbert polynomials. \par

\smallskip

\noindent\textbf{Acknowledgements.} We thank Ciro Ciliberto, Patricio Gallardo, Rick Miranda, Anand Patel and Debaditya Raychaudhury for motivating discussions and for generously sharing time listening to the results. The first author also thanks Rick Miranda and Jeanne Duflot for arranging his visit and hospitality at Fort Collins. We also thank Rick Miranda and Helge Ruddat for pointing out some relevant references and for their helpful comments that improved the exposition.

\color{black}

\smallskip

\section{Main Results}\label{main results}

We assume throughout that we work over an algebraically closed field of characteristic zero. In this section we develop the methods that are the foundations for the rest of the paper. We start with the following lemma that relates the Hilbert polynomial of certain embedded simple normal crossings to the Hilbert polynomial of an embedded ribbon.

\begin{lemma}\label{Hilbert polynomial}
Let $Y$ be a smooth projective scheme embedded (possibly degenerately) inside a projective space 
 $\mathbb{P}^N$ and let $D \in \mid L^{-1} \mid$ be a smooth divisor in $Y$. Suppose there exist
\begin{itemize}
    \item[(1)] A scheme $V$ inside $\mathbb{P}^N$ which is the (scheme-theoretic) union $Y_1 \bigcup Y_2$ of two smooth schemes inside $\mathbb{P}^N$ such that both $Y_1$ and $Y_2$  belong to the same Hilbert scheme as $Y$ and $Y_1 \bigcap Y_2 = D$ (scheme-theoretic intersection) as a subscheme of $\mathbb{P}^N$. Then,  
    $$ \chi(\mathscr{O}_V(t)) = \chi(\mathcal{O}_{Y_1}(t)) +  \chi(\mathcal{O}_{Y_2}(t)) - \chi(\mathcal{O}_{D}(t)) $$
      
     and, for $t>>0$,
     $$ \chi(\mathscr{O}_V(t)) = 2\chi(\mathcal{O}_{Y}(t))-\chi(\mathcal{O}_{D}(t)).$$
    \item[(2)] A ribbon $\widetilde{Y}$ on $Y$ with conormal bundle $L$ embedded inside $\mathbb{P}^N$. Then, for $t>>0$, $$ \chi(\mathscr{O}_{\widetilde{Y}}(t)) = 2\chi(\mathcal{O}_{Y}(t))-\chi(\mathcal{O}_{D}(t)).$$
\end{itemize}
In particular, the Hilbert polynomials of $V$ and $\widetilde{Y}$ are the same. 
\end{lemma}

\noindent\textbf{Proof.} Let $\mathscr{O}_Y(1)$, $\mathscr{O}_{Y_i}(1)$ and $\mathscr{O}_D(1)$ be the line bundles obtained by pulling back $\mathscr{O}_{\mathbb{P}^N}(1)$ to the respective varieties. Let $S = \mathbf{k}[x_0,..,x_N]$ and $I(Y_i)$ denote the homogeneous ideal of $Y_i$.It follows from the, well-known exact sequence

$$ 0 \to \displaystyle\frac{S}{I(Y_1) \bigcap I(Y_2)} \to \displaystyle\frac{S}{I(Y_1)} \bigoplus \displaystyle\frac{S}{I(Y_2)} \to \displaystyle\frac{S}{(I(Y_1) + I(Y_2))} \to 0 $$

we get

\color{black}

$$ \chi(\mathscr{O}_V(t)) = \chi(\mathcal{O}_{Y_1}(t)) +  \chi(\mathcal{O}_{Y_2}(t)) - \chi(\mathcal{O}_{D}(t)) .$$
Since both $Y_1$ and $Y_2$ belong to the same Hilbert scheme as $Y$, $Y_1$, $Y_2$ and $Y$ have the same Hilbert polynomial, so, for $t>>0$, we have
$$\chi(\mathcal{O}_{Y_1}(t)) +  \chi(\mathcal{O}_{Y_2}(t)) - \chi(\mathcal{O}_{D}(t)) = 2\chi(\mathcal{O}_{Y}(t))-\chi(\mathcal{O}_{D}(t)).$$
This proves (1).\\
Now, from the exact sequence
$${0 \to L \to \mathcal{O}_Y \to \mathcal{O}_D \to 0}.$$

we get $\chi(\mathcal{O}_{Y}(t)) - \chi(\mathcal{O}_{D}(t)) = \chi(L \otimes \mathcal{O}_{Y}(t))$.
Then, it follows from (1) that, for $t >>0$ $$ \chi(\mathcal{O}_{V}(t)) = \chi(\mathcal{O}_{Y}(t))+ \chi(L \otimes \mathcal{O}_{Y}(t)).$$
Now we have the exact sequence
\begin{equation}\label{sequence.conormal.ribbon}
0 \to L \to \mathcal{O}_{\widetilde{Y}} \to \mathcal{O}_Y \to 0 
\end{equation}
which yields
$$ \chi(\mathcal{O}_{\widetilde{Y}}(t)) =  \chi(\mathcal{O}_{Y}(t))+ \chi(L \otimes \mathcal{O}_{Y}(t)).$$ This proves (2)\QEDA

\color{black}

\begin{remark}\label{remark.Hartshorne}(\cite{H10}, Proposition $2.3$)\label{ideal of deformation}
 Let $Y \hookrightarrow \mathbb{P}^N$ be an embedding of a smooth variety inside a projective space. Let $\mathcal{I}$ denote the ideal sheaf of $Y$. Let $Y_{\Delta}$ denote an embedded first order deformation of $Y$ given by $\Bar{\psi} \in \textrm{Hom}(\mathcal{I}/\mathcal{I}^2, \mathcal{O}_Y) = H^0(N_{Y/\mathbb{P}^N})$.
  Then the ideal $\mathcal{I}_{\Delta}$ of $Y_{\Delta}$ inside $\mathbb{P}_{\Delta}^N$ is given over an affine open $\textrm{Spec}(R) \times \Delta \subset \mathbb{P}_{\Delta}^N $ (where $\textrm{Spec}(R) \subset \mathbb{P}^N$ and $\Delta = k[\epsilon]/\epsilon^2$) by
 $$ I_{\Delta} = \{a+b\epsilon \mid a \in I, b \in R, 
  \psi(a) = \Bar{b}  \} \subset R \times k[\epsilon]/\epsilon^2 $$ where $I \subset R$ is the ideal of $Y$ restricted to affine open $\textrm{Spec}(R)$ and $\Bar{b}$ is the image of $b$ under the canonical quotient $R \to R/I$ and $\psi$ is the composition $\psi: I \to I/I^2 \xrightarrow{\Bar{\psi}} R/I$. \\
  In particular, if $Y_{\Delta} = Y \times \Delta$ is the trivial first order deformation of $X$ corresponding to $\Bar{\psi} = 0 \in \textrm{Hom}(\mathcal{I}/\mathcal{I}^2, \mathcal{O}_Y)$, then the ideal of $Y_{\Delta}$ is given by 
  $$ I_{\Delta} = \{a+b\epsilon \mid a \in I, b \in I 
\} \subset R \times k[\epsilon]/\epsilon^2$$
or $$ I_{\Delta} = I \bigoplus I\epsilon \subset R \times k[\epsilon]/\epsilon^2 $$
\end{remark}

\begin{remark}\label{ideal of a ribbon}
(see the proof of \cite[Lemma 1.4]{GP97}; see also \cite{HV})
Let $Y \hookrightarrow \mathbb{P}^N$ be an embedding of a smooth variety inside a projective space. Let $\mathcal{I}$ denote the ideal sheaf of $Y$. Let $\widetilde{Y}$ denote an embedded ribbon on $Y$ with conormal bundle $L$ corresponding to a surjective homomorphism $\Bar{\psi} \in \textrm{Hom}(\mathcal{I}/\mathcal{I}^2, L) = H^0(N_{Y/\mathbb{P}^N} \otimes L)$. Then the ideal sheaf $\mathcal{I}_{\widetilde{Y}}$ of $\widetilde{Y}$ inside $\mathbb{P}^N$ is given by $\textrm{Ker}\psi$ where $\psi$ is the composition of $\psi: \mathcal{I} \to \mathcal{I}/\mathcal{I}^2 \to L$
\end{remark}

\begin{lemma}\label{deformations fixing a subscheme}
Let $D \hookrightarrow Y \hookrightarrow \mathbb{P}^N$ be the embedding of smooth subvariety $D$ and variety  $Y$ inside a projective space over an algebraically closed field $\mathbf k$ of characteristic zero. Let $\textrm{Art}_\mathbf k$ denote the category of artinian $\mathbf k$-algebras. Let $F: \textrm{Art}_\textbf k \to \textrm{Sets}$ be the functor as follows: For an Artinian $\mathbf k-$algebra $A$,
$$F_{D,Y}(A) = \left\{ D \times \textrm{Spec}(A) \hookrightarrow Y_A \hookrightarrow \mathbb{P}^N \times \textrm{Spec}(A) \right\}$$
where $Y_A \hookrightarrow \mathbb{P}^N \times \textrm{Spec}(A)$ is an embedded deformation of $Y$ inside $\mathbb{P}^N$ over $\textrm{Spec}(A)$. Then 
\begin{itemize}
    \item[(1)] $F_{D,Y}(\mathbf k[\epsilon]/\epsilon^2) = H^0(N_{Y/\mathbb{P}^N} \otimes \mathcal{I}_{D/Y})$ where $\mathcal{I}_{D/Y}$ is the ideal sheaf of $D$ inside $Y$.
    \item[(2)] $H^1(N_{Y/\mathbb{P}^N} \otimes \mathcal{I}_{D/Y})$ is an obstruction space for the functor $F_{D,Y}$. 
\end{itemize}
\end{lemma}

\noindent\textit{Proof.} The functor $F_{D,Y}$ parametrizing subschemes $Y$ of $\mathbb P^N$ that contain $D$, is a natural generalization of the Hilbert scheme (see \cite[p. 230]{Ser}). Then the result follows as a natural generalization of \cite[Theorem VI-29]{EH} or the analogous result in \cite[Proposition 3.2.6]{Ser} or \cite{H10}. \QEDA

\begin{remark}\label{flag hilbert schemes}
Let $i: Y \hookrightarrow \mathbb{P}^N$ be an embedding of a smooth variety inside a projective space. Let $D \in |L^{-1}|$ be a smooth divisor. Then there is an induced embedding of $j: D \hookrightarrow \mathbb{P}^N$. Let $p_1$ and $p_2$ denote the Hilbert polynomials of $Y \hookrightarrow \mathbb{P}^N$ and $D \hookrightarrow \mathbb{P}^N$ respectively. Let $\mathcal{H}(p_1, p_2)$ denote the flag Hilbert scheme representing the flag-Hilbert functor as defined in \cite{Ser}, Section $4.5.1$ or \cite{KL81}. Closed points of $\mathcal{H}(p_1, p_2)$ parameterize pairs $(D' \hookrightarrow Y' \hookrightarrow \mathbb{P}^N)$ where the subschemes $Y' \hookrightarrow \mathbb{P}^N$ and $D' \hookrightarrow \mathbb{P}^N$ have Hilbert polynomials $p_1$ and $p_2$ respectively. Let $\mathcal{D}$ denote the Hilbert scheme parameterizing subschemes $(D' \hookrightarrow \mathbb{P}^N)$ with Hilbert polynomial $p_2$. There is a projection map $p: \mathcal{H}(p_1, p_2) \to \mathcal{D}$ whose fibre $\mathcal{F}_D$ at a subscheme $(D \hookrightarrow \mathbb{P}^N)$ consists of all those subschemes $Y' \hookrightarrow \mathbb{P}^N$ with Hilbert polynomial $p_1$ such that $D \subset Y'$. It follows from the definition of flag-Hilbert functor $\mathcal{H}(p_1,p_2)$ given in \cite{Ser} or \cite{KL81}, that $F_{D,Y}$ defined in Lemma \ref{deformations fixing a subscheme} is represented by the pointed scheme consisting of  $(\mathcal{F}_D, (D \hookrightarrow Y \hookrightarrow \mathbb{P}^N))$. In particular the dimension of the fibre at $(D \hookrightarrow Y \hookrightarrow \mathbb{P}^N)$ is $H^0(N_{Y/\mathbb{P}^N} \otimes L)$ and the fibre is smooth at the point $(D \hookrightarrow Y \hookrightarrow \mathbb{P}^N)$ if $H^1(N_{Y/\mathbb{P}^N} \otimes L) = 0$. Let $\mathcal F_{D,Y}$ denote the unique irreducible component of $\mathcal F_D$ containing the point $(D \hookrightarrow Y\hookrightarrow \mathbb{P}^N)$. Closed points parameterized by $\mathcal{F}_{D,Y}$ correspond to deformations over an irreducible curve of the subvariety $Y$ inside $\mathbb{P}^N$ which keeps $D$ fixed.     
\end{remark}

\color{black}

We now prove the first of the two important theorems of this section. In this theorem, we use the existence of embedded ribbons to construct a one-parameter family of simple normal crossing varieties that degenerate to the ribbon. \newline
Let $Y$ be a smooth projective variety embedded inside 
$\mathbb{P}^N$ and let $D \in |L^{-1}|$ be a smooth divisor. Let $p_1$ and $p_2$ denote the Hilbert polynomials $Y$ and $D$ inside $\mathbb{P}^N$ respectively. Let $\mathcal{D}, p: \mathcal{H}(p_1,p_2) \to \mathcal{D}$ and $F_{D,Y}$ be as in Lemma \ref{deformations fixing a subscheme} and Remark \ref{flag hilbert schemes}. \\

\smallskip

\begin{definition}
    Let $Y_1$ and $Y_2$ be two smooth subvarieties of $\mathbb P^N$ and let $D$ be a smooth divisor on both $Y_1$ and $Y_2$. If $D$ is the scheme-theoretic intersection of $Y_1$ and $Y_2$, we denote the union $Y_1 \bigcup Y_2$ as $Y_1 \bigcup_D Y_2$.
\end{definition}

\color{black}
\begin{theorem}\label{main}
Let $Y$ be a smooth projective variety embedded in $\mathbb{P}^N$. Let $L$ be a line bundle on $Y$ and  
$D \in \mid L^{-1} \mid$  be a smooth divisor. Suppose there exists a nowhere vanishing section $\zeta$ in $H^0(N_{Y/\mathbb{P}^N} \otimes L)$ which corresponds to 
a ribbon 
$\widetilde{Y}$, embedded in $\mathbb{P}^N$, supported on $Y$. Assume 
that $Y$ is an unobstructed subscheme of $\mathbb{P}^N$, and that the functor $F_{D,Y}$ is unobstructed or equivalently the fibre of the map $p: \mathcal{H}(p_1,p_2) \to \mathcal{D}$ over the point $(D \hookrightarrow \mathbb{P}^N)$ is smooth at the point $(D \hookrightarrow Y \hookrightarrow \mathbb{P}^N)$. Then there exists a flat family of subschemes $\mathcal{Y} \to T$ of $\mathbb{P}^N_T$  over a smooth irreducible algebraic curve $T$ such that 
\begin{itemize}
    \item[(1)]  $\mathcal{Y}_0$ is the subscheme $\widetilde{Y}$. 
    \item[(2)] $\mathcal{Y}_t = Y_{1t} \bigcup Y_{2t}$, where $Y_{1t}$ is the subscheme $Y$, $Y_{2t}$ is a non trivial deformation in $\mathbb{P}^N_T$ of 
    the subscheme $Y$ that fixes the subscheme $D$
    (i.e., 
    for any $t \in T, t \neq 0$, $Y_{2t}$ and $Y$ are different subschemes of $\mathbb{P}^N_T$)
    and $Y_{1t} \bigcap Y_{2t} \cong D$.
\end{itemize}
\end{theorem}

\noindent\textbf{Proof.} \underline{Step $1$: Construction of the family:} \textcolor{black}{Let $\lambda$ be the sheaf inclusion 
$$\lambda: L \to \mathcal O_Y$$
induced by $D \in \mid L^{-1} \mid$. Making an abuse of notation, we} consider the exact sequence
$$0 \to H^0(N_{Y/\mathbb{P}^N}\otimes L) \xrightarrow{\lambda} H^0(N_{Y/\mathbb{P}^N}).$$
Now consider \textcolor{black}{the subscheme $\lambda(\widetilde{Y})$ of $\mathbb P^N_\Delta$ corresponding to $\lambda(\zeta)$.}

By Lemma \ref{deformations fixing a subscheme}, \textcolor{black}{The subscheme $\lambda(\widetilde{Y})$} corresponds to an embedded first order deformation that contains the trivial embedded first order deformation of $D \hookrightarrow \mathbb{P}^N$. \\ Since  \textcolor{black}{the functor $F_{D,Y}$ is unobstructed or equivalently the fibre of the map $p: \mathcal{H}(p_1,p_2) \to \mathcal{D}$ over the point $(D \hookrightarrow \mathbb{P}^N)$ is smooth at the point $(D \hookrightarrow Y \hookrightarrow \mathbb{P}^N)$}, $\lambda(\widetilde{Y})$  can be extended to a deformation $\mathcal{Y}_2 \hookrightarrow \mathbb{P}^N_T$ (smooth and irreducible) over a smooth curve $T$ containing the trivial deformation of $D \hookrightarrow \mathbb{P}^N$, i.e.,
$$ D \times T \hookrightarrow \mathcal{Y}_2 \hookrightarrow \mathbb{P}^N_T.$$
\textcolor{black}{Note that, since $\zeta$ is non zero, $\mathcal{Y}_2$ is non trivial in the sense of the statement.}\\

Now consider the trivial deformation (smooth and irreducible)  $\mathcal{Y}_1 = Y \times T \hookrightarrow \mathbb{P}_T^N$. Note that  we  have once again
$$ D \times T \hookrightarrow \mathcal{Y}_1 \hookrightarrow \mathbb{P}_T^N  $$   \\
Consider $\mathcal{Y}_1 \bigcup \mathcal{Y}_2 \to T$. We claim that this family satisfies the conditions of our theorem. \\
 
\underline{Step 2. Flatness of the family:} Both irreducible components $\mathcal{Y}_1$ and $\mathcal{Y}_2$ of $\mathcal{Y}_1 \bigcup \mathcal{Y}_2$ dominate $T$. \textcolor{black}{Since both $\mathcal{Y}_1$ and $\mathcal{Y}_2$ are reduced, so is $\mathcal{Y}_1 \bigcup \mathcal{Y}_2$. Thus $\mathcal{Y}_1 \bigcup \mathcal{Y}_2$ is flat over $T$.} 

\underline{Step 3. The central fibre \textcolor{black}{of the family} contains the ribbon $\widetilde{Y}$:} Let us look at the first order deformation induced by $\mathcal{Y}_1 \bigcup \mathcal{Y}_2$ which is $(\mathcal{Y}_1 \bigcup \mathcal{Y}_2)_{\Delta} \supset \mathcal{Y}_{1\Delta} \bigcup \mathcal{Y}_{2\Delta} = (Y \times \Delta) \bigcup \lambda({\widetilde{Y}})$. Now recall from Remark ~\ref{ideal of deformation} (or Step $1$) that
$$I_{\lambda{(\widetilde{Y})}} = \{a+b\epsilon \mid a \in I, b \in R, {(\lambda \circ \psi)}(a) = \Bar{b} \} \subset R \times k[\epsilon]/\epsilon^2 $$ 
Also
$$ I_{Y \times \Delta} = \{a+b\epsilon \mid a \in I, b \in I \} \subset R \times k[\epsilon]/\epsilon^2 $$
We have $$I_{\lambda({\widetilde{Y}}) \bigcup (Y \times \Delta)} = I_{\lambda({\widetilde{Y}})} \bigcap I_{(Y \times \Delta)}$$ \\
Claim:  $J = \{a+b\epsilon \mid a \in \textrm{Ker}{(\lambda \circ \psi)}, b \in I \} = I_{\lambda({\widetilde{Y}})} \bigcap I_{(Y \times \Delta)} \subset R \times k[\epsilon]/\epsilon^2 $. \\
\noindent\textit{Proof:} Note that if $a+b\epsilon \in J$, we have ${(\lambda \circ \psi)}(a) = 0$ and $\Bar{b} = 0$ and hence $J \subset I_{\lambda{(\widetilde{Y})}}$. Again since both $a, b \in I$ (by definition), we have $J \subset I_{Y \times \Delta}$. So $J \subseteq I_{\lambda({\widetilde{Y}})} \bigcap I_{(Y \times \Delta)}$. On the other hand, suppose $a+b\epsilon \in I_{\lambda({\widetilde{Y}})} \bigcap I_{(Y \times \Delta)}$. Then $a \in I, b \in I$ and ${(\lambda \circ \psi)}(a) = \Bar{b} $. But $ b \in I$, hence $\Bar{b} = 0$ and ${(\lambda \circ \psi)}(a) = 0$, so  $a \in \textrm{Ker}{(\lambda \circ \psi)}$. This implies $J \supseteq I_{\lambda({\widetilde{Y}})} \bigcap I_{(Y \times \Delta)}$.   \\
Hence the central fibre of the family contains 
$$ (\lambda({\widetilde{Y}}) \bigcup (Y \times \Delta))_0 $$ whose ideal is given inside $\textrm{Spec}(R)$ by $J/\epsilon J = \textrm{Ker}(\lambda \circ \psi) = \textrm{Ker}(\psi)$ (since $\lambda$ is injective). But $\textrm{Ker}(\psi)$ is by Remark ~\ref{ideal of a ribbon} the ideal of $\widetilde{Y}$. \\ 

\underline{Step 4. The central fibre is $\widetilde{Y}$ and the general fibre is a variety with normal crossing:} \\
For a general $t \in T$,
\textcolor{black}{the fiber of the family  is 
$\mathcal{Y}_{1t} \bigcup \mathcal{Y}_{2t}$, with $\mathcal{Y}_{1t}$ different from $\mathcal{Y}_{2t}$ and}  $\mathcal{Y}_{1t} \bigcap \mathcal{Y}_{2t}$ containing the divisor $D$ as a subscheme. \\
\color{black}
Suppose $\mathcal{Y}_{1t} \bigcap \mathcal{Y}_{2t} = D'$. Then 
$$0 \to \mathcal{J} \to \mathcal{O}_{D'} \to \mathcal{O}_D \to 0$$ where $\mathcal{J}$ is the ideal sheaf of $D$ inside $D'$. Then $$\chi(\mathcal{O}_{D'}(t)) = \chi(\mathcal{O}_{D}(t)) + \chi(\mathcal{J}(t)).$$ Then 
\begin{equation}
\chi(\mathcal{O}_{\mathcal{Y}_{1t} \bigcup \mathcal{Y}_{2t}}(t)) = \chi(\mathcal{O}_{\mathcal{Y}_{1t}}(t)) + \chi(\mathcal{O}_{\mathcal{Y}_{2t}}(t)) - \chi(\mathcal{O}_{\mathcal{Y}_{1t} \bigcap \mathcal{Y}_{2t}}(t) = \chi(\mathcal{O}_{\mathcal{Y}_{1t}}(t)) + \chi(\mathcal{O}_{\mathcal{Y}_{2t}}(t)) - \chi(\mathcal{O}_{D'}(t))
\end{equation}
and hence, \textcolor{black}{for $t>>0$,}
\begin{equation}\label{one1}
  \chi(\mathcal{O}_{\mathcal{Y}_{1t} \bigcup \mathcal{Y}_{2t}}(t)) = 2\chi(\mathcal{O}_{Y}(t)) - \chi(\mathcal{O}_{D}(t)) - \chi(\textcolor{black}{\mathcal{J}}(t))
\end{equation}
Now \textcolor{black}{recall} that, \textcolor{black}{by Step 3}, the central fibre is a scheme $X$ containing $\widetilde{Y}$ as a subscheme, with ideal sheaf $\mathcal{K}$. Then we have 
$$0 \to \mathcal{K} \to \mathcal{O}_X \to \mathcal{O}_{\widetilde{Y}} \to 0$$ and hence, using Lemma ~\ref{Hilbert polynomial}, \textcolor{black}{for $t>>0$,}
\begin{equation}\label{two2}
\chi(\mathcal{O}_{X}(t)) = \chi(\mathcal{O}_{\widetilde{Y}}(t)) + \chi(\mathcal{K}(t)) = 2\chi(\mathcal{O}_{Y}(t)) - \chi(\mathcal{O}_{D}(t)) + \chi(\mathcal{K}(t))
\end{equation}
Now note that since $X$ and $\mathcal{Y}_{1t} \bigcup \mathcal{Y}_{2t}$ are fibres of a flat family their Hilbert polynomials are the same and hence by equations ~\ref{one1} and ~\ref{two2}, both $\chi(\mathcal{K}(t))$ and $\chi(\textcolor{black}{\mathcal{J}}(t))$ are equal to $0$ for $t >> 0$ which implies that $X = \widetilde{Y}$ and $\mathcal{Y}_{1t} \bigcap \mathcal{Y}_{2t} = D$.
Now since $\mathcal{Y}_{1t}$ and $ \mathcal{Y}_{2t}$ are both smooth and they intersect along the smooth divisor $D$ for a general $t$, it follows that $\mathcal{Y}_{1t} \bigcup_D \mathcal{Y}_{2t}$ is simple normal crossing (by \href{https://stacks.math.columbia.edu/tag/0CBN}{Tag 0CBN} ,Lemma $41.21.2$, $(2)$). \QEDA 

\color{black}

\begin{remark}\label{blow-up of a ribbon}
    Let $\widetilde{Y}$ be a ribbon over $Y$ with conormal bundle $L$. Let $D \in |L^{-1}|$ denote a smooth divisor on $Y$. Then $D$ defines a Weil-divisor on $\widetilde{Y}$. Choosing a section $s$ corresponding to $D$, we have a map $H^1(T_Y \otimes L) \to H^1(T_Y)$ obtained by multiplication by the section $s$. Note that $\widetilde{Y}$ represents an element of $H^1(T_Y \otimes L)$ and under the above map the image of $\widetilde{Y}$ is the blow-up of $\widetilde{Y}$ along the Weil divisor $D$, which is a ribbon with conormal bundle $L \otimes D = \mathcal{O}_Y$ (see \cite[Theorem 1.9]{BE95}) . A ribbon with conormal bundle $\mathcal{O}_Y$ is by the definition of a ribbon (see Definition \ref{def of a ribbon}) a first order deformation of $Y$, so indeed the image lies $H^1(T_Y)$. Now in the embedded setting, we have a commutative diagram
    \[
    \begin{tikzcd}
        H^1(N_{Y/\mathbb{P}^N} \otimes L) \arrow[d] \arrow[r, "\lambda"] & H^1(N_{Y/\mathbb{P}^N}) \arrow[d] \\
        H^1(T_Y \otimes L) \arrow[r] & H^1(T_Y) \\
    \end{tikzcd}
    \]
where the vertical maps are the forgetful maps. So for the embedded ribbon $\widetilde{Y} \hookrightarrow \mathbb{P}^N$, the image $\lambda(\widetilde{Y})$ is abstractly the blow-up of $\widetilde{Y}$ along $D$ and hence is itself an embedded (in $\mathbb{P}_{\Delta}^N$) first order deformation of $Y$.
\end{remark}

\begin{definition}
    Let $Y$ be a smooth projective variety. A smooth projective variety $Y'$ is deformation equivalent to $Y$ if there exists an irreducible curve $T$ and a flat family $\mathcal{Y}$ such that $\mathcal{Y}_t = Y$ and $\mathcal{Y}_{t'} = Y'$ for two different points $t$ and $t'$. 
\end{definition}

In the rest of the article, we will follow the conventions and notations set up in Lemma \ref{deformations fixing a subscheme} and Remark \ref{flag hilbert schemes}. The following is the second main theorem of this section where we combine Theorem ~\ref{main} with the techniques on smoothing of ribbons developed in \cite{Go} and \cite{GGP13}, to give an algorithm on how to construct simple normal crossing varieties and carry out their smoothings inside a projective space. We also compute the codimension of the locus of the equisingular deformations of such smoothable simple normal crossings. 

\begin{theorem}\label{intrinsic}
    Let $Y$ be a smooth projective variety and $D \in |L^{-1}|$ be a smooth divisor. Assume there exists an embedding $i: Y \hookrightarrow \mathbb{P}^N$ such {that the} following conditions hold{:}
    \begin{itemize}
        \item[(a)] There exists a  {nowhere vanishing section} {$\nu$} in $H^0(N_{Y/\mathbb{P}^N} \otimes L)$. 
        \item[(b)] The variety {$Y$ is unobstructed in $\mathbb P^N$.}
        \item[(c)] The functor $F_{D,Y}$ is unobstructed  
         (note that this holds if $H^1(N_{Y/\mathbb{P}^N} \otimes L) = 0$). 
        \item[(d)] {The embedded ribbon}  $\widetilde{Y}$ {corresponding to $\nu$} is smoothable in $\mathbb{P}^N$.
        \item[(e)] The Hilbert point of $\widetilde{Y}$ in $\mathbb{P}^N$ is smooth. 
    \end{itemize}
Then:
\begin{itemize}
    \item[(1)] 
   If $Y_1$ is the subvariety $Y$ of $\mathbb P^N$ and $Y_2$ is a subvariety of $\mathbb P^N$ that corresponds to a general element of 
    $\mathcal{F}_{D,Y}$, then {the scheme theoretic intersection of $Y_1 and Y_2$ is $D$ and} 
{$V = Y_1 \bigcup_D Y_2$} is  smoothable inside $\mathbb{P}^N$. Moreover, the Hilbert point of $V$ inside $\mathbb{P}^N$ is smooth. 
\end{itemize}

{Let $N_{V/\mathbb{P}^N}'$ denote the equisingular 
normal sheaf of $V$ inside $\mathbb{P}^N$ and assume that, in addition to conditions (a), (b), (c), (d) and (e),  the cohomology group $H^1(N_{V/\mathbb{P}^N}')$ vanishes. Then:} 
\begin{itemize}
    \item[(2)] The locally trivial deformations of $V$ inside $\mathbb{P}^N$ are unobstructed
    and the first order locally trivial deformations form a subspace in the tangent space of the Hilbert scheme at $[V]$ of codimension $h^0(L^{-2}|_D)$. 
    Furthermore, this subspace is the tangent space of an  irreducible locus of the Hilbert scheme, which has codimension $h^0(L^{-2}|_D)$ and is smooth at $[V]$. The singularities of the subschemes parameterized by this locus are  normal crossing singularities;  {in particular, they are non-normal and semi log canonical. They are analytically isomorphic to $(x_1^2+x_2^2 = 0 \subset \mathbb{C}^N)$.}
\end{itemize}

\end{theorem}

\noindent\textit{Proof:} 
{First  we prove} {the existence of a nontrivial deformation $Y_2$ of $Y$ in $\mathbb P^N$ satisfying the thesis of (1).} Since the conditions of Theorem \ref{main} are fulfilled, we have a flat family of subschemes $\mathcal{Y} \to T$ of $\mathbb{P}^N_T$ over a smooth irreducible algebraic curve $T$ such that 
\begin{itemize}
    \item[(i)] $\mathcal{Y}_0$ is the subscheme $\widetilde{Y}$. \item[(ii)] $\mathcal{Y}_t = Y_{1t} \bigcup Y_{2t}$, for each $t \neq 0$, where $Y_{1t}$ is the subscheme $Y$, {and the varieties $Y_{2t}$}  
 {form} a non trivial deformation in $\mathbb{P}^N$ of the subscheme $Y$ and {the intersection of $Y_{1t}$  and $Y_{2t}$ is the subscheme $D$.} 
\end{itemize} 
Since $T$ is irreducible, there exists an irreducible component of the Hilbert scheme of $\mathbb{P}^N$ containing both $\widetilde{Y}$ and $\mathcal{Y}_t = Y_{1t} \bigcup Y_{2t}$ for each $t \neq 0$. Now by condition $(e)$, $\widetilde{Y}$ is contained in a unique component of the Hilbert scheme and by $(d)$, a general element of that component is a smooth subscheme of $\mathbb{P}^N$. \newline 
Now we prove {that the subvariety $Y_2$ above can be chosen to be a general element of $\mathcal{F}_{D,Y}$.}  
Let $\mathfrak Y \longrightarrow \mathcal F_{D,Y}$ be the universal family over $\mathcal F_{D,Y}$. Consider 
$$\left(Y \times \mathcal F_{D,Y}\right) \bigcup \mathfrak Y \overset{\Phi}\longrightarrow \mathcal F_{D,Y}.$$
There is a non-empty open set $U$  of $\mathcal F_{D,Y}$ {over which $\Phi$ is flat and the fibers over the closed points of $U$ are the form  
$Y \bigcup Y_2'$, where {$Y_2'$ corresponds to a point of $\mathcal F_{D,Y}$.} 
Maybe after shrinking, the family  $\mathcal Y \longrightarrow T$ in the proof of (1) may be taken so that $T \setminus \{0\}$ is contained in $U$ ({and the fiber over $t \in T, t \neq 0$ is the same as the fiber over $t$ of the family parameterized by $U$}). Then the fibers of $\Phi$ over the points of $U$ have the same Hilbert polynomial as the ribbon $\widetilde Y$, so the fibers of $\Phi$ are normal crossing varieties $Y \bigcup_D Y_2'$, where the varieties $Y_2'$ are parameterized by $U$. Since $\mathcal F_{D,Y}$ is irreducible, arguing as in the proof of (1) we conclude that all the varieties  $Y \bigcup_D Y_2'$, with $Y_2'$ parameterized by $U$, can be smoothed. Since the Hilbert point of $\widetilde Y$ is smooth, after shrinking $U$ if necessary, we may conclude that the Hilbert point of such $Y \bigcup_D Y_2'$ is also smooth. \color{black}

Now we prove (2). 
{Since $H^1(N_{V/\mathbb{P}^N}')$ is an obstruction space for the locally trivial deformations of $V$ in $\mathbb{P}^N$ (see \cite[Proposition 1.1.9, Theorem 3.4.8, Corollary 3.4.9]{Ser}), {the locally trivial deformations of $V$ are unobstructed}. Now recall that $H^0(N_{V/\mathbb{P}^N}')$ is the tangent space to the locally trivial deformations of $V$ in $\mathbb P^N$  and that $N_{V/\mathbb{P}^N}'$ fits into the exact sequence
$$ 0    \longrightarrow N_{V/\mathbb{P}^N}' \longrightarrow N_{V/\mathbb{P}^N} \longrightarrow T_V^1 \longrightarrow 0$$
(see \cite[\S 4.7.1]{Ser}), so we have 
$$0 \to H^0(N_{V/\mathbb{P}^N}') \to H^0(N_{V/\mathbb{P}^N}) \to H^0(T_V^1) \to 0.$$
This implies {our claim about codimension inside the tangent to the Hilbert scheme}, 
since $T_V^1 = L^{-2}|_D$ (see \cite[Proposition 2.3]{F83}). 
Then the locus described in (2) is obtained from the algebraization of the formal semiuniversal locally trivial deformation of $V$ in $\mathbb P^N$ \textcolor{black}{(see \cite[Theorem 3.4.8]{Ser})} and the codimension of this locus and its smoothness at $[V]$ follow from the above. Finally, at its singular points, $V$ is locally analytically isomorphic to the zero locus of $x_1^2+x_2^2$ in $\mathbb C^N$, and, by \cite[Corollary 11.2]{TN10}, so are the subvarieties corresponding to general points of the locus. Then, in particular, they are normal crossing, non-normal and semi log canonical subvarieties.} \QEDA

\begin{theorem}\label{theorem.irreducible}
Let $Y$ and $D$ satisfy all the hypotheses of Theorem~\ref{intrinsic}. Assume furthermore that hypotheses (a) and (d) hold for a general element of  $\mathcal{F}_{D,Y}$ ({note this holds if $H^1(N_{Y/\mathbb{P}^N} \otimes L) = 0$, see Remark~\ref{remark.generality}}). {Then for a general pair $(Y_1',Y_2')$ in  $\mathcal{F}_{D,Y} \times \mathcal{F}_{D,Y}$, the union of $Y_1'$ and $Y_2'$ is a normal crossing subvariety {$V=Y_1' \bigcup_D Y_2'$} of $\mathbb{P}^N$, which}
    is smoothable. Such subvarieties $V$ form  an irreducible, flat family. In particular, if $\mathcal{F}_{D}$ is irreducible  {and hypotheses (a) and (d) are fulfilled by a general element of  $\mathcal{F}_{D}$, then the above holds for a general pair $(Y_1',Y_2')$ in  $\mathcal{F}_{D} \times \mathcal{F}_{D}$}.
\end{theorem}
    
\noindent\textit{Proof.} Let $\mathfrak Y \longrightarrow \mathcal F_{Y,D}$ the restriction to $\mathcal F_{Y,D}$ of the universal family parameterizing subschemes containing $D$. {Consider the family $$\left[\left[(\mathfrak{Y} \times \mathbb{P}_{\mathcal F_{Y,D}}^N) \bigcup (\mathbb{P}_{\mathcal F_{Y,D}}^N \times \mathfrak{Y})\right] \bigcap (\Delta \times \mathcal F_{Y,D} \times \mathcal F_{Y,D})\right]_{\textrm{red}}$$ over $\mathcal F_{Y,D} \times \mathcal F_{Y,D}$.} Let $\mathcal U$ be an open set of 
$\mathcal F_{Y,D} \times \mathcal F_{Y,D}$ on which the universal family is  flat {and reduced} and {hence} the scheme theoretic fiber over any closed point $(s,t) \in \mathcal U, s \neq t$, is $\mathfrak Y_s \bigcup \mathfrak Y_t$. By hypothesis, there exists an open  set {of $\mathcal F_{Y,D}$}  where hypotheses (a) and (d)  hold. We may assume that (b), (c) and (e) also hold on that open set because smoothness is an open condition. Let  $\mathcal U'$ be the intersection of the open set of $\mathcal F_{Y,D}$ where hypotheses (a) to (e) hold with $p_1(\mathcal U)$, where $p_1$ is the projection from $\mathcal F_{Y,D} \times \mathcal F_{Y,D}$ onto its first factor. {We want to show that {the varieties} $\mathfrak Y_{s} \bigcup \mathfrak Y_{t}$, where $(s,t) \in p_1^{-1}\mathcal U' \bigcap \mathcal U, s \neq t$, are normal crossing varieties $\mathfrak Y_{s} \bigcup_D \mathfrak Y_{t}$ and are smoothable. To see this,}
fix $(s_0,t_0) \in p_1^{-1}\mathcal U' \bigcap \mathcal U, s_0 \neq t_0$ and rename $\mathcal Y_{s_0}$ as $Y$ and  let $$U= p_1^{-1}\mathcal U' \bigcap \mathcal U \bigcap \left(\{s_0\} \times \mathcal F_{Y,D}\right).$$ Arguing as in the proof of Theorem~\ref{intrinsic} (2), we conclude that the varieties $\mathcal Y_{s_0} \bigcup \mathcal Y_{t}$, where $t \in U, t \neq 0$ have the same Hilbert polynomial as $\widetilde Y$, {and hence}  so do the varieties $\mathcal Y_{s} \bigcup \mathcal Y_{t}$, where $(s,t) \in p_1^{-1}\mathcal U' \bigcap \mathcal U, s \neq t$, {since $p_1^{-1}\mathcal U' \bigcap \mathcal U$ is irreducible}. Thus such 
varieties are normal crossing varieties $\mathcal Y_{s} \bigcup_D \mathcal Y_{t}$ and are smoothable. \QEDA

\begin{remark}\label{remark.generality}
Let $\mathcal{F}_{D}$ be irreducible
\begin{enumerate}
    \item A sufficient condition for hypotheses (a) and (d) to hold for a general element of $\mathcal{F}_{D}$ is $h^0(N_{\mathcal Y_t/\mathbb P^N} \otimes  L)$ to be constant for general elements of $\mathcal{F}_{D}$. 
    \item A sufficient condition for $h^0(N_{\mathcal Y_t/\mathbb P^N} \otimes  L)$ to be constant for general elements of $\mathcal{F}_{D}$ is $h^1(N_{\mathcal Y_t/\mathbb P^N} \otimes  L)$ to be $0$. 
\end{enumerate}
\end{remark}

\color{black}

\smallskip

{In the following remark, we give sufficient conditions for condition (a) of Theorem ~\ref{intrinsic} to hold.}

\begin{remark}(\textbf{A criterion for the existence of 
nowhere vanishing global sections of $N_{Y/\mathbb{P}^N} \otimes L$)}.\label{existence}
    Let $Y$ be a smooth projective variety of dimension $d$ and $D \in |L^{-1}|$ be a smooth divisor in $Y$. Let $H$ be a very ample line bundle and consider an embedding $Y \xhookrightarrow{i} \mathbb{P}^N$ (possibly degenerate) such that $i^*(\mathcal{O}_{\mathbb{P}^N}(1)) = H$. Suppose
    \begin{itemize}
        \item[(1)] $H \otimes L$ is base point free 
        \item[(2)] $h^0(H) \geq 2d+2$
    \end{itemize}
    Then there exists nowhere vanishing sections in $H^0(N_{Y/\mathbb{P}^N} \otimes L)$. 
\end{remark}

\noindent\textit{Proof.}  Consider the Euler sequence of $\mathbb{P}^N$ restricted to $Y$ and twisted by $L$.
$$ 0 \to L \to (H \otimes L)^{N+1} \to T_{\mathbb{P}^N}|_Y \otimes L \to 0 $$
Since $H \otimes L$ is base point free and since $T_{\mathbb{P}^N}|_Y \otimes L$ is non zero (it has rank $N$) we conclude that $T_{\mathbb{P}^N}|_Y \otimes L$ is base point free.    
Now consider the exact sequence 
$$0 \to T_Y \otimes L \to T_{\mathbb{P}^N}|_Y \otimes L \to N_{Y/\mathbb{P}^N} \otimes L \to 0$$ 
Now $N_{Y/\mathbb{P}^N} \otimes L$ is non-zero and has rank $N-d > d$, it follows that $N_{Y/\mathbb{P}^N} \otimes L$ is base point free and has a nowhere vanishing section.
\smallskip

{The following lemma addresses the condition (e) of Theorem ~\ref{intrinsic}. It \textcolor{black}{deduces} the condition of smoothness of Hilbert point of the ribbons \textcolor{black}{from} a series of cohomology vanishing conditions on its reduced part.}

\begin{lemma}\label{smoothness of Hilbert point}\textbf{(A criterion for the smoothness of the Hilbert point of a ribbon)}
Let $Y \hookrightarrow \mathbb{P}^N$ be an embedding where $Y$ is smooth. Let $L$ be a line bundle on $Y$ and let $\widetilde{Y} \hookrightarrow \mathbb{P}^N$ be an embedding of a ribbon $\widetilde{Y}$ on $Y$ with conormal bundle $L$ extending the embedding $Y \hookrightarrow \mathbb{P}^N$. Suppose the following conditions hold
\begin{itemize}
    \item[(1)] $H^1(L^{-2}) = H^1(N_{Y/\mathbb{P}^N}) = H^2(L^*) = 0$
    \item[(2)] $H^1(L^{*}) = H^1(N_{Y/\mathbb{P}^N}\otimes L) = H^2(\mathcal{O}_Y) = 0$
    
\end{itemize}
Then $H^1(N_{\widetilde{Y}/\mathbb{P}^N}) = 0$ and consequently the Hilbert point of $\widetilde{Y}$ inside $\mathbb{P}^N$ is smooth.
\end{lemma}

\noindent\textit{Proof.} We have the following exact sequences
\begin{equation}\label{1}
    0 \to N_{\widetilde{Y}/\mathbb{P}^N} \otimes L \to N_{\widetilde{Y}/\mathbb{P}^N} \to N_{\widetilde{Y}/\mathbb{P}^N} \otimes \mathcal{O}_Y \to 0
\end{equation}
\begin{equation}\label{2}
    0 \to \textrm{Hom}(I_{\widetilde{Y}}/I_{Y}^2, \mathcal{O}_Y) \to N_{\widetilde{Y}/\mathbb{P}^N} \otimes \mathcal{O}_Y \to L^{-2} \to 0 
\end{equation}
\begin{equation}\label{3}
    0 \to L^* \to N_{Y/\mathbb{P}^N} \to \textrm{Hom}(I_{\widetilde{Y}}/I_{Y}^2, \mathcal{O}_Y) \to 0
\end{equation}
Tensoring ~\ref{2} by $L$ we have 
\begin{equation}\label{4}
    0 \to \textrm{Hom}(I_{\widetilde{Y}}/I_{Y}^2, \mathcal{O}_Y) \otimes L \to N_{\widetilde{Y}/\mathbb{P}^N} \otimes L \to L^{*} \to 0
\end{equation}
Tensoring ~\ref{3} by $L$ we have 
\begin{equation}\label{5}
    0 \to \mathcal{O}_Y \to N_{Y/\mathbb{P}^N} \otimes L \to \textrm{Hom}(I_{\widetilde{Y}}/I_{Y}^2, \mathcal{O}_Y) \otimes L \to 0
\end{equation}
Using exact sequence ~\ref{2} and ~\ref{3} and $H^1(L^{-2}) = H^1(N_{Y/\mathbb{P}^N}) = H^2(L^*) = 0$ we have that $H^1(N_{\widetilde{Y}/\mathbb{P}^N} \otimes \mathcal{O}_Y) = 0$.
Using exact sequence ~\ref{4} and ~\ref{5} and $H^1(L^{*}) = H^1(N_{Y/\mathbb{P}^N}\otimes L) = H^2(\mathcal{O}_Y) = 0$ we have that $H^1(N_{\widetilde{Y}/\mathbb{P}^N} \otimes L) = 0$. Now using ~\ref{1} we conclude that $H^1(N_{\widetilde{Y}/\mathbb{P}^N}) = 0$. By \cite{Ser}, Proposition $3.2.6$, since $\widetilde{Y}$ is a local complete intersection, we have that the Hilbert point of $\widetilde{Y}$ in $\mathbb{P}^N$ is smooth. \QEDA

\begin{remark}\label{cohomology of normal}
Let $Y \xhookrightarrow{i} \mathbb{P}^N$ be a possibly degenerate embedding of a smooth projective variety $Y$ such that $i^*(\mathcal{O}_{\mathbb{P}^N}(1)) = \mathcal{O}_Y(1)$. Let $L$ be a line bundle on $Y$.
\begin{itemize}
    \item[(1)] Assume $H^2(\mathcal{O}_Y) = H^1(\mathcal{O}_Y(1)) = H^2(T_Y) = 0$. Then $H^1(N_{Y/\mathbb{P}^N}) = 0$.   
    
    \item[(2)]  Assume $H^2(L) = H^1(\mathcal{O}_Y(1) \otimes L) = H^2(T_Y \otimes L) = 0$. Then $H^1(N_{Y/\mathbb{P}^N} \otimes L) = 0$.

\end{itemize}
    
\end{remark}

\noindent\textit{Proof.} Consider the restriction of the Euler sequence of $\mathbb{P}^N$ to $Y$ 
$$0 \to \mathcal{O}_Y \to \mathcal{O}_Y(1)^{N+1} \to T_{\mathbb{P}^N}|_Y \to 0$$ Note that since $H^2(\mathcal{O}_Y) = H^3(\mathcal{O}_Y) = H^1(\mathcal{O}_Y(1)) = H^2(\mathcal{O}_Y(1)) = 0$, $H^1(T_{\mathbb{P}^N}|_Y) = H^2(T_{\mathbb{P}^N}|_Y) = 0$. Now part $(1)$ follows when we consider the exact sequence
$$0 \to T_Y \to T_{\mathbb{P}^N}|_Y \to N_{Y/\mathbb{P}^N} \to 0$$
Consider the restriction of the Euler sequence of $\mathbb{P}^N$ to $Y$ twisted by $L$. 
$$0 \to L \to (\mathcal{O}_Y(1) \otimes L)^{N+1} \to T_{\mathbb{P}^N}|_Y \otimes L \to 0$$ Note that since $H^2(L) = H^3(L) = H^1(\mathcal{O}_Y(1) \otimes L) = H^2(\mathcal{O}_Y(1) \otimes L) = 0$, $H^1(T_{\mathbb{P}^N}|_Y \otimes L) = H^2(T_{\mathbb{P}^N}|_Y \otimes L) = 0$. Now part $(2)$ follows when we consider the exact sequence
$$0 \to T_Y \otimes L \to T_{\mathbb{P}^N}|_Y \otimes L \to N_{Y/\mathbb{P}^N} \otimes L \to 0$$

{Now we address the condition (d) of Theorem \ref{intrinsic}, by formulating Proposition $1.4$ of \cite{GGP13} in a form which is convenient for us.}

\begin{proposition}(see \cite{GGP13}) \label{thm.criterion.smoothing} \textbf{(A criterion for embedded smoothing of ribbons)}\label{smoothing}
 Let $Y \hookrightarrow \mathbb{P}^N$ be an embedding where $Y$ is smooth. Let $L$ be a line bundle on $Y$ and let $\widetilde{Y} \hookrightarrow \mathbb{P}^N$ be an embedding of a ribbon $\widetilde{Y}$ on $Y$ with conormal bundle $L$ extending the embedding $Y \hookrightarrow \mathbb{P}^N$. Suppose the following conditions hold
\begin{itemize}
    \item[(1)] There exists a smooth divisor $B \in |L^{-2}|$.
    \item[(2)] $H^1(B|_B) = 0$ {(this happens e.g. if  $H^1(B) = H^2(\mathcal{O}_Y) = 0$).} 
    \item[(3)] $H^1(N_{Y/\mathbb{P}^N}) = H^1(N_{Y/\mathbb{P}^N}\otimes L) = 0$
    
\end{itemize}
Then $\widetilde{Y}$ is smoothable inside $\mathbb{P}^N$.    
\end{proposition}

\noindent\textit{Proof.}   Since $-2L$ has a smooth member $B$, we can consider the smooth double cover $X \xrightarrow{\pi} Y$ branched along the smooth divisor $B$. Now consider the composition $\varphi = i \circ \pi$ where $i: Y \hookrightarrow \mathbb{P}^N$ is the embedding. Now consider the exact sequence 
$$ 0 \to N_{\pi} \to N_{\varphi} \to \pi^*(N_{Y/\mathbb{P}^N}) \to 0 $$
After taking the cohomology we have the exact sequence $$H^0(N_{\varphi}) \to H^0(N_{Y/\mathbb{P}^N}) \oplus H^0(N_{Y/\mathbb{P}^N}\otimes L) \to H^1(N_{\pi})$$
Note that $H^1(N_{\pi}) = H^1(B|_B) = 0$ by the conditions of the theorem. Hence the induced map $H^0(N_{\varphi}) \to H^0(N_{Y/\mathbb{P}^N}\otimes L)$ surjects. Now once again we have an exact sequence
$$H^1(N_{\pi}) \to H^1(N_{\varphi}) \to H^1(N_{Y/\mathbb{P}^N}) \oplus H^1(N_{Y/\mathbb{P}^N}\otimes L) \to 0$$
Again by the conditions of the theorem, $H^1(N_{\varphi}) = 0$ and hence the morphism $\varphi$ is unobstructed. Therefore by \cite{GGP13}, Proposition $1.4$, $\widetilde{Y}$ is smoothable. \QEDA

\smallskip

Theorem ~\ref{intrinsic}, gives us sufficient conditions to show existence and smoothing of simple normal crossings inside a projective space. In the next proposition, we use Remark ~\ref{existence}, Lemma ~\ref{smoothness of Hilbert point}, Proposition ~\ref{smoothing} to reduce those conditions into checking the vanishing of the cohomology and positivity of certain line bundles on the components. 
In the proposition we use the notation introduced in Remark 2.3 for the flag Hilbert Scheme.

\begin{proposition}\label{abstract smoothing 1}
    Let $Y$ be a smooth projective variety. Let $D \in |L^{-1}|$ be a smooth divisor. Assume the following conditions are satisfied:
    \begin{itemize}
        \item[(1)] $H^1(L^{-2}) = H^2(L^*) = H^1(L^{*}) = H^2(\mathcal{O}_Y) = H^2(L) = H^2(T_Y) = 0$. 
        \item[(2)] The homomorphism 
        $H^2(T_Y \otimes L) \to H^2(T_{\mathbb P^N}|_Y \otimes L)$ arising from the long exact sequence of cohomology of the normal exact sequence is injective.
        \item[(3)] There exists a smooth divisor $B \in |L^{-2}|$.
    \end{itemize}
Suppose there exists a very ample line bundle $H$ on $Y$ inducing an embedding $Y \xhookrightarrow{i} \mathbb{P}^N$, such that
\begin{itemize}
        \item[(i)] $H \otimes L$ is base point free 
        \item[(ii)] $h^0(H) = N+1 \geq 2d+2$
        \item[(iii)] $H^1(H \otimes L) = 0$ 
        \item[(iv)] $H^1(H) = 0$ 
    \end{itemize}   

If $Y_1$ is the subvariety $Y$ of $\mathbb P^N$ and $Y_2$ is a subvariety of $\mathbb P^N$ that corresponds to a general element of $\mathcal{F}_{D,Y}$, then {the scheme theoretic intersection of $Y_1 and Y_2$ is $D$ and} {$V = Y_1 \bigcup_D Y_2$} is  smoothable inside $\mathbb{P}^N$. Moreover, the Hilbert point of $V$ inside $\mathbb{P}^N$ is smooth.

\end{proposition}

\noindent\textit{Proof.} We check the conditions of Theorem ~\ref{intrinsic}. By Remark \ref{existence}, if we consider the embedding $Y \hookrightarrow \mathbb{P}^N$ by the complete linear series of the very ample line bundle $H$, then condition $(1)$ of Theorem ~\ref{intrinsic} is satisfied and there exist embedded ribbons $\widetilde{Y}$ inside $\mathbb{P}^N$ which extends the embedding $Y \hookrightarrow \mathbb{P}^N$. Now since $H^1(H \otimes L) = H^2(L) = 0$ \textcolor{black}{and by condition (2),}  we conclude that $H^1(N_{Y/\mathbb{P}^N}\otimes L) = 0$. To see that $\widetilde{Y}$ is smoothable we check the conditions of Theorem ~\ref{smoothing}.  First of all, note that there exists a smooth divisor $B \in |-2L|$. Under the assumptions, we have that $H^1(B) = H^1(-2L) = H^2(\mathcal{O}_Y) = 0$. We have already shown that $H^1(N_{Y/\mathbb{P}^N}\otimes L) = 0$. Again under the assumptions, $H^2(\mathcal{O}_Y) = H^1(H) = H^2(T_Y) = 0$. Hence by Remark ~\ref{cohomology of normal}, $(1)$, we have that $H^1(N_{Y/\mathbb{P}^N}) = 0$. Hence by Theorem ~\ref{smoothing}, $\widetilde{Y}$ is smoothable. Now we check the smoothness of the Hilbert point $\widetilde{Y}$ inside $\mathbb{P}^N$ using Lemma ~\ref{smoothness of Hilbert point}. Note that under the assumptions, $H^1(L^{-2}) = H^2(L^*) = 0$ and we just showed that $H^1(N_{Y/\mathbb{P}^N}\otimes L) = 0$. Similarly under the conditions, $H^1(L^{*}) = H^2(\mathcal{O}_Y) = 0$ and we also showed that $H^1(N_{Y/\mathbb{P}^N}\otimes L) = 0$. Hence by Lemma ~\ref{smoothness of Hilbert point}, the Hilbert point of $\widetilde{Y}$ inside $\mathbb{P}^N$ is smooth. Therefore the result follows from Theorem ~\ref{intrinsic}. \QEDA

\section{Projective smoothing of union of Fano varieties}\label{Fano of any dimension}

In this section we apply the general {theory and} methods  {developed in Section~\ref{main results}} to show the existence of {simple normal crossing subvarieties $V$ of $\mathbb{P}^N$ which are the union of two irreducible Fano varieties $Y_1$ and $Y_2$.} {We show that these subvarieties $V$ can be smoothed inside $\mathbb P^N$.} We start with the following definition:

\begin{definition}\label{Kodaira dimension of normal crossing}
   \begin{itemize}
    \item[(1)] A variety $X$ is Fano if it is semi-log canonical, $X$ is projective and $-K_X$ is an ample $\mathbb{Q}-$ divisor.
    \item[(2)]  A variety $X$ is called Calabi-Yau, if it is semi-log canonical, $X$ is proper, $K_X$ is $\mathbb{Q}-$ linearly equivalent to $0$ and $H^i(\mathcal{O}_X) = 0$ for $0 < i < \textrm{dim}(X)$. 
    \item[(3)] A variety $X$ is stable if it is semi-log canonical, $X$ is projective and $K_X$ is an ample $\mathbb{Q}-$ divisor. 
\end{itemize}  
\end{definition}

\begin{theorem}\label{Fano}
{Let $Y$ be a smooth, projective Fano variety of dimension $d \geq 3$, let $H$ be a very ample line bundle on $Y$, let  $N=h^0(H)-1$ and, abusing the notation, call $Y$ the image in $\mathbb P^N$ of the embedding induced by the complete linear series $|H|$. Let $L$ be a line bundle such that $L^{-1}$ is ample and there exist a smooth divisor  $D$ in  $|L^{-1}|$ and a smooth divisor  in $|L^{-2}|$. Assume the following conditions hold: 
\begin{enumerate}
     \item[(a)] $H \otimes L$ is base point free.
     \item[(b)] $N \geq 2d+1$. 
     \item[(c)] One of the following conditions hold:
     \begin{enumerate}
        \item[(i)] $\omega_Y^{-1} \otimes L$ is ample; or
        \item[(ii)] $L =  \omega_Y$ and $h^{1,1}(Y) = 1$ when $d = 3$ or
        $L= \omega_Y$ and $h^{1,d-2}(Y) = 0$
       when $d \geq 4$; 
        \item[(iii)]  $L^{-1} \otimes  \omega_Y $ is ample, 
        $H^1(H  \otimes L)=0$
        and $Y$ satisfies Bott vanishing. 
    \end{enumerate}
\end{enumerate}
Then:  
\begin{itemize}
    \item[(1)] There exist simple normal crossing subvarieties $V = Y_1 \bigcup_D Y_2$ {in $\mathbb{P}^N$} as in Theorem~\ref{intrinsic} (1);  {in particular}, $V$ can be {smoothed} inside $\mathbb{P}^N$. 
    \item[(2)] More precisely, the subvarieties $V$ are Fano varieties, Calabi-Yau varieties or stable varieties and can be deformed inside $\mathbb{P}^N$ to smooth Fano varieties, Calabi-Yau varieties or varieties with ample canonical bundle {in cases (i), (ii), (iii) respectively}. 
    \item[(3)] The locally trivial deformations of $V$ inside $\mathbb{P}^N$  are as in Theorem~\ref{intrinsic} (2).
    \item[(4)] The general fiber of any one-parameter family of deformations of $V$ whose image under the Kodaira-Spencer map is a first order locally trivial deformation of $V$, is once again a simple normal crossing variety $V' = Y_1' \bigcup_D' Y_2'$, where both $Y_1'$ and $Y_2'$ are deformations of the subvariety $Y$ and $D'$ is a deformation of the subvariety $D$. 
\end{itemize}}
\end{theorem}

\color{black}

\vspace{0.2cm}

\noindent\textit{Proof.} {We will use Proposition ~\ref{abstract smoothing 1} to prove (1), so we check that the hypotheses of Proposition ~\ref{abstract smoothing 1} hold. The existence of smooth divisors in $|L^{-1}|$ and $|L^{-2}|$, the fact that $H \otimes L$ is base point free and the bound $N \geq 2d+1$ are satisfied by hypothesis. For the remaining conditions of Proposition ~\ref{abstract smoothing 1}, recall that  $d \geq 3$ and $\omega_Y^{-1}$, $L^{-1}$ and $H$ are ample. Then the vanishings of $H^1(L^{-1})$,  $H^2(L^{-1})$, 
$H^1(L^{-2})$, $H^2(\mathcal O_Y)$, $H^2(L)$ and $H^1(H)$ follow from the Kodaira vanishing theorem, and the vanishing of $H^2(T_Y)$ follows from the Nakano vanishing theorem. 
The vanishing of  $H^1(H \otimes L)$ follows from Kodaira vanishing taking into account the assumptions in (i), (ii) and (iii). \\
\smallskip
We now check the only hypothesis of Proposition~\ref{abstract smoothing 1} left, namely, condition (2). Assume first that  $d=3$ and $L=\omega_Y$ does not occur. Then  
$H^2(T_Y \otimes L)$ vanishes, by the Nakano vanishing theorem in case (i), by Serre duality in case (ii) and from Bott vanishing in case (iii), so Condition (2) of Proposition~\ref{abstract smoothing 1} is satisfied.
\\
Now assume $d=3$ and $L=\omega_Y$. Note that, in this case, $H^2(T_Y \otimes L) \neq 0$, because $h^{1,1}(Y) \neq 0$, so we cannot argue as above. We then look more closely at the homomorphism of Proposition~\ref{abstract smoothing 1} (2), which in our situation is
\begin{equation*}
  H^2(T_Y \otimes \omega_Y) \overset{f} {\longrightarrow} H^2(T_{\mathbb{P}^N}|_Y \otimes \omega_Y).
\end{equation*}
The homomorphism $f$ can be identified, via Serre duality, with the dual of the homomorphism
\begin{equation}\label{$(1)$}
    H^1(\Omega_{\mathbb{P}^N}|_Y) \overset{g} {\longrightarrow} H^1(\Omega_Y).
\end{equation}
Composing $g$ with the natural map 
$H^1(\Omega_{\mathbb{P}^N}) \to H^1(\Omega_{\mathbb{P}^N}|_Y)$, we get 
a homomorphism 
\begin{equation*}
    H^1(\Omega_{\mathbb{P}^N}) \overset{\widehat g} {\longrightarrow} H^1(\Omega_Y).
\end{equation*}
which induces, by restriction,  a homomorphism between the Neron-Severi groups of $\mathbb P^N$ and $Y$, which, in our situation, is a homomorphism $h$ between the Picard groups of $\mathbb P^N$ and $Y$, because $Y$, being Fano, is a regular variety. \textcolor{black}{The homomorphism $h$ sends a line bundle on $\mathbb P^N$ to its restriction on $Y$}, so it is nonzero. Then, so is $\widehat g$ and, consequently, so is $g$. Therefore, $f$ is also nonzero. Since, by Serre duality, $h^2(T_Y \otimes \omega_Y) = h^1(\Omega_Y)$ and, by assumption, $h^{1,1}(Y) = 1$, the homomorphism $f$ is injective as desired (in fact, it is an isomorphism). Then Condition (2) of Proposition 2.17 is satisfied in all cases, so Proposition 2.17 implies Theorem~\ref{Fano} (1).}  

\smallskip
{\textcolor{black}{Now we} prove (2). 
\color{black}
It is enough to prove that $V$ is Fano, Calabi-Yau or stable in cases (i), (ii) and (iii) respectively. Note that $\omega_V|_{Y_i} = \omega_{Y_i}(D)$ and $\omega_V$ is negative ample or trivial or ample if and only if so is $\omega_V|_{Y_i}$. Now $\omega_{Y_i}(D)$ is negative ample under condition (i), trivial under condition (ii) and ample under condition (iii). The fact that under condition (ii) all intermediate cohomologies of $\mathcal{O}_V$ vanish, follow from the exact sequence since $Y$ is Fano and $L^{-1}$ is ample  
$$ 0 \to L \to \mathcal{O}_V \to \mathcal{O}_{Y_2} \to 0 $$

\smallskip

{In order to prove (3), notice that, by Theorem~\ref{intrinsic}, it suffices to show that $H^1(N_{V/\mathbb{P}^N}') = 0$.
Recall the exact sequence
$$0 \to T_V \to T_{\mathbb{P}^N}|_V \to N_{V/\mathbb{P}^N}' \to 0$$ (see \cite[Proposition 1.1.9]{Ser}), 
that yields the exact sequence
\begin{equation}\label{cohomology.normal.sequence.V}
 H^1(T_{\mathbb{P}^N}|_V) \to H^1(N_{V/\mathbb{P}^N}') \to H^2(T_V) \overset{f_V}{\longrightarrow} H^2(T_{\mathbb{P}^N}|_V). 
\end{equation}
First we see that $H^1(T_{\mathbb{P}^N}|_V)$ vanishes.  
Since $V=Y_1 \cup Y_2$ and $D$ is the scheme theoretic intersection on $Y_1$ and $Y_2$, we have the following exact sequence of $\mathcal O_V$-modules: 
\begin{equation}\label{equation.O_V}
0 \to L \to \mathcal{O}_V \to \mathcal{O}_{Y_2} \to 0  
\end{equation}
In the proof of (1) we have already seen that $H^1(H \otimes L)=H^2(L)$. We have also seen that $H^1(H)=H^2(\mathcal O_Y)=0$. Since $Y_2$ is a general member of a flat, one parameter family of deformations of $Y$, this implies that 
$H^1(\mathcal O_{Y_2}(1))$ and $H^2(\mathcal O_{Y_2})$ vanish.
Then twisting and taking cohomology on 
\eqref{equation.O_V} we get $H^1(\mathcal O_V(1))=H^2(\mathcal O_V)=0$.
Then taking cohomology on the Euler sequence of $\mathbb P^N$ restricted to $V$ yields $H^1(T_{\mathbb{P}^N}|_V)=0$. 
\\
\smallskip
Now we look at the homomorphism $f_V$ of \eqref{cohomology.normal.sequence.V}. Except if $d=3$ and $L=\omega_Y$, $H^2(T_V)$ vanishes by Lemma~\ref{H^2(T_V)} below. Now assume $d=3$ and $L=\omega_Y$. By \cite [Corollary 4.2]{TN15} (see also \cite[Lemma 2.9]{F83}), we have $$H^2(T_V) = H^{1}(\Omega_V/\tau_V \otimes \omega_V)^*,$$ where $\tau_V$ is the torsion subsheaf of $\Omega_V$. By Serre duality,  the homomorphism $f_V$  of \ref{cohomology.normal.sequence.V} is the dual of the map 
\begin{equation}\label{equation.cohomology.torsion.free}
    H^1(\Omega_{\mathbb{P}^N}|_{V} \otimes \omega_V) \longrightarrow  H^1((\Omega_V/\tau_V) \otimes \omega_V).
\end{equation} 
Since, by Theorem~\ref{intrinsic}, $V$  is a general member of a flat, one parameter family of deformations of a ribbon $\widetilde Y$ supported on $Y$, with conormal bundle $L^{-1} = \omega_Y$. Then $\omega_{\widetilde Y}$ is  trivial (see the proof of \cite[Corollary 4.4]{BGG}), so $\omega_V$ is also trivial. 
Then the map in \eqref{equation.cohomology.torsion.free} is  
$$H^1(\Omega_{\mathbb{P}^N}|_{V}) \overset{g_V}{\longrightarrow} H^1(\Omega_V/\tau_V).$$ From the exact sequence \cite[(4.1)]{TN15}, we focus on the map
\begin{equation}\label{equation.Tzi}
 \Omega_V/\tau_V \to (\pi_0)_*(\Omega_{\widetilde{V}}),  
\end{equation}
where $\pi_0: \widetilde{V} \to V$ is the normalization of $V$, so $\widetilde{V}$ is the disjoint union of $Y_1$ and $Y_2$. Taking cohomology on \eqref{equation.Tzi}, we get  the map
$$H^{1}(\Omega_V/\tau_V)  \overset{g_V'}{\longrightarrow} H^1((\pi_0)_*(\Omega_{\widetilde{V}})) \simeq \bigoplus_{j = 1}^2 H^{1}(\Omega_{Y_j}).$$
Composing $g_V$ with the natural map  $H^1(\Omega_{\mathbb{P}^N}) \to H^1(\Omega_{\mathbb{P}^N}|_{V})$ on the left and the map  $g'_V$ on the right we obtain a map 
$$H^1(\Omega_{\mathbb{P}^N}) \longrightarrow \bigoplus_{j = 1}^2 H^{1}(\Omega_{Y_j}),$$  which 
induces a map 
$$NS(\mathbb{P}^N) {\longrightarrow} \bigoplus_{j = 1}^2 NS(Y_j),$$ which, since $Y_1$ and $Y_2$ are regular varieties, is in fact a map 
$$\mathrm{Pic}(\mathbb{P}^N) \overset{h_V}{\longrightarrow} \bigoplus_{j = 1}^2 \mathrm{Pic}(Y_j).$$
\textcolor{black}{The homomorphism $h_V$ sends a line bundle $\mathcal L$ on 
$\mathbb P^N$ to $\pi_0^*(\mathcal L|_V)$}, so 
$h_V$ is nonzero. Therefore $g_V$ is nonzero and so is $f_V$. Since by Lemma~\ref{H^2(T_V)}, 
$h^2(T_V)=1$, then $h_V$ is injective )in fact, it is an isomorphism). This, together with the vanishing of $H^1(T_{\mathbb P^N}|_V)$, which we have already proved, implies $H^1(N_{V/\mathbb P^N}')$ vanishes, so (3) is proved.} 

To \textcolor{black}{prove (4)}, note that given a simple normal crossing variety $V = Y_1 \bigcup_D Y_2$  \textcolor{black}{as} in part $(1)$ {(recall \textcolor{black}{that, in particular,} $Y_1 = Y$)}, there is a sublocus of embedded deformations of $V$, whose Kodaira-Spencer map lands in the subspace of locally trivial deformations, which consists of simple normal crossing varieties $V' = Y_1' \bigcup_D' Y_2'$, where both $Y_1'$ and $Y_2'$ are deformations of the subscheme $Y$ and $D'$ is a deformation of the subscheme $D$. We calculate the dimension of this  \textcolor{black}{sub}locus and show that it is the same as the dimension of the \textcolor{black}{locally trivial} locus. One can calculate the dimension by the following steps: 
\begin{itemize}
    \item[(a)] Choose a general $Y_1'$ from the irreducible component of the Hilbert scheme containing $Y = Y_1$.
    \item[(b)] For each such $Y_1'$ choose a smooth divisor $D'$ from the linear system of $L'$, which is the unique deformation of $L$ on $Y$ to $Y_1'$ since $H^2(\mathcal{O}_Y) = H^1(\mathcal{O}_Y) = 0$.  
    \item[(c)] Now one can apply Theorem \ref{intrinsic} $(1)$ to the pair $(Y_1', D')$ to conclude that a general $Y_2'$ in $\mathcal{F}_{D', Y_1'}$ intersects $Y_1'$ at exactly $D'$. 
\end{itemize}
This gives a smooth irreducible family of simple normal crossing varieties $V' = Y_1' \bigcup_D' Y_2'$ containing $V = Y \bigcup_D Y_2$, where both $Y_1'$ and $Y_2'$ are deformations of the subscheme $Y$ and $D'$ is a deformation of the subscheme $D$ of dimension of $h^0(N_{Y/\mathbb{P}^N})+h^0(L^{-1})-1+h^0(N_{Y/\mathbb{P}^N} \otimes L)$.
On the other hand the dimension of the ambient Hilbert component, which is the Hilbert component containing the ribbon $\widetilde{Y}$ is $h^0(N_{\widetilde{Y}/\mathbb{P}^N})$, which by Lemma \ref{smoothness of Hilbert point}, can be shown to be $h^0(N_{Y/\mathbb{P}^N})+h^0(N_{Y/\mathbb{P}^N} \otimes L)+h^0(L^{-2})-1$. So the codimension is given by $h^0(L^{-2})-h^0(L^{-1}) = h^0(L^{-2}|_D)$ which is the codimension of the locally trivial locus as shown before in Theorem \ref{intrinsic} $(2)$.  
\QEDA

\medskip

\color{black}

In the following Lemma we compute $H^2(T_V)$ of the normal crossing varieties $V$ that we construct in Theorem \ref{Fano}. This was required for the proof of Theorem \ref{Fano}, part $(3)$.

\begin{lemma}\label{H^2(T_V)}
In the situation of Theorem \ref{Fano},
\begin{itemize}
    \item[(1)] \textcolor{black}{if $d=3$ and $L=\omega_Y$, then $h^2(T_V) = 1$;}
 
    \item[(2)] \textcolor{black}{otherwise, $h^2(T_V) = 0$.} 
\end{itemize}
    
\end{lemma}

\noindent\textbf{Proof.} In Case \textcolor{black}{(i) of Theorem~\ref{Fano}}, the result follows from \textcolor{black}{Theorem~\ref{Fano} (2) and \cite[Theorem 4.8]{TN15}}, Theorem $4.8$}. For the other cases, first note that from \cite[Corollary 4.2]{TN15}  (see also \cite[Lemma 2.9]{F83}), we have 
$$H^2(T_V) = H^{n-2}(\Omega_V/\tau_V \otimes \omega_V)^*.$$ 
\color{black}
Recall that $Y$ is a smooth Fano variety of dimension $d \geq 3$ and $D$ is a smooth, ample subvariety of $Y$. By \cite[Corollary 4.18]{Debarre} and \cite[Theorem 0.1]{KMM},  $Y$ is simply connected, so, by the Lefschetz theorem  for homotopy groups, so is $D$. Then $D$ does not have any torsion line bundles other than $\mathcal O_D$. Therefore, after tensoring with $\omega_V$, the exact sequence  \cite[(4.4.1)]{TN15} becomes
\begin{equation}\label{modulo.torsion.sequence}
  0 \to \Omega_V/\tau_V \otimes \omega_V \to (\pi_0)_*(\Omega_{\widetilde{V}}) \otimes \omega_V \to \Omega_{D} \otimes \omega_V \to 0,
\end{equation} 
where $\pi_0: \widetilde{V} \longrightarrow V$ is the normalization of $V$.
Since $D$ is the scheme--theoretic intersection of $Y_1$ and $Y_2$, we have 
\begin{equation}\label{canonical.restricted.to.Y}
 \omega_V|_{Y_i} = \omega_{Y_i}(D)    
\end{equation}
and hence \begin{equation}\label{canonical.restricted.to.D}
    \omega_V|_D = 
 \omega_V(D)|_D = \omega_D.
\end{equation}
 \color{black}
On the other hand by projection formula and because $\pi_0$ is finite and $\widetilde V$ is isomorphic to the disjoint union of $Y_1$ and $Y_2$, we have for any $l \geq 0$ 
\begin{equation}\label{cohomology.equation}
    H^{l}((\pi_0)_*(\Omega_{\widetilde{V}}) \otimes \omega_V) = H^{l}((\pi_0)_*(\Omega_{\widetilde{V}} \otimes \pi_0^*\omega_V)) = H^{l}(\Omega_{\widetilde{V}} \otimes \pi_0^*\omega_V) = \bigoplus_{j = 1}^2 H^{l}(\Omega_{Y_j} \otimes \omega_{Y_j}(D)).
\end{equation}
By \eqref{canonical.restricted.to.Y}, \eqref{canonical.restricted.to.D} and \eqref{cohomology.equation}, 
after taking cohomology on \eqref{modulo.torsion.sequence}
we get 
\begin{multline*}
    \bigoplus_{j = 1}^2 H^{n-3}(\Omega_{Y_j} \otimes \omega_{Y_j}(D)) \to H^{n-3}(\Omega_{D} \otimes  \omega_D) \to  H^{n-2}(\Omega_V/\tau_V \otimes \omega_V) \\ \to \bigoplus_{j = 1}^2 H^{n-2}(\Omega_{Y_j} \otimes \omega_{Y_j}(D)) \to H^{n-2}(\Omega_{D} \otimes \omega_D). 
\end{multline*}
\color{black}

First we show that $$H^{n-3}((\Omega_{Y_j} \otimes (\omega_{Y_j}(D)))) \to H^{n-3}(\Omega_{D} \otimes \omega_D)$$ is surjective. Since $Y_j \cong Y$, we have to show that the map $$H^{n-3}((\Omega_{Y} \otimes (\omega_{Y}(D)))) \to H^{n-3}(\Omega_{D} \otimes \omega_D)$$ is surjective. 
We have the conormal exact sequence
$$0 \to -D|_D \to \Omega_Y|_D \to \Omega_D \to 0$$
Tensoring with $\omega_{Y}(D)$ we get 
$$0 \to \omega_Y|_D \to \Omega_Y|_D \otimes (\omega_{Y}(D)) \to \Omega_D \otimes \omega_D \to 0 $$
Now consider the exact sequence 
$$0 \to -D \to \mathcal{O}_Y \to \mathcal{O}_D \to 0 $$
Tensoring the sequence by $\Omega_Y \otimes (\omega_{Y}(D))$ we have 
$$0 \to \Omega_Y \otimes \omega_Y \to \Omega_Y \otimes \omega_Y(D) \to \Omega_Y|_D \otimes \omega_Y(D)  \to 0 $$
Since the map $\Omega_{Y} \otimes \omega_Y(D) \to \Omega_{D} \otimes \omega_D$ factors through $\Omega_Y|_D \otimes \omega_Y(D)$ we have from the above two exact sequence that the map $\Omega_{Y} \otimes (\omega_{Y}(D)) \to \Omega_{D} \otimes \omega_D$ is surjective. 
Letting $\mathcal{K}$ being the kernel of $\Omega_{Y} \otimes (\omega_{Y}(D)) \to \Omega_{D} \otimes \omega_D$ we have a commutative diagram
\[
\begin{tikzcd}
    & & 0 \arrow[d] & & \\
    & & \Omega_Y \arrow[d] \otimes \omega_Y & & \\
    0 \arrow[r]  & \mathcal{K} \arrow[r] \arrow[d, "f"] & \Omega_Y \otimes \omega_Y(D) \arrow[r] \arrow[d, "g"] & \Omega_{D} \otimes \omega_D \arrow[r] \arrow[d, "\cong"] & 0 \\
    0 \arrow[r]  & \omega_Y|_D \arrow[r] & \Omega_Y|_D \otimes \omega_Y(D) \arrow[r] \arrow[d] & \Omega_{D} \otimes \omega_D \arrow[r] \arrow[d] & 0 \\
    & & 0 & 0
\end{tikzcd}
\]
Applying snake lemma, we have that $\textrm{Ker}(f) = \Omega_Y \otimes \omega_Y$ and $\textrm{Coker}(f) = 0$. Hence we have an exact sequence
$$ 0 \to \Omega_Y \otimes \omega_Y \to \mathcal{K} \to \omega_Y|_D \to 0 $$
Now note that show the surjection of $H^{n-3}((\Omega_{Y} \otimes (\omega_{Y}(D)))) \to H^{n-3}(\Omega_{D} \otimes \omega_D)$ we need to show that 
$H^{n-2}(\mathcal{K}) = 0$. Considering the last exact sequence it is enough to show $H^{n-2}(\Omega_Y \otimes \omega_Y) = H^{n-2}(\omega_Y|_D) = 0$. Since $Y$ is Fano, the first vanishing follows by Nakano vanishing theorem. For the second one, consider the exact sequence
$$0 \to \omega_Y(-D) \to \omega_Y \to \omega_Y|_D \to 0$$
Note that again since $Y$ is Fano, $H^{n-2}(\omega_Y) = 0$. In case (ii) and (iii), we have $\omega_Y(-D) = \omega_Y \otimes L$ which is negative ample and hence $H^{n-1}(\omega_Y(-D)) = 0$. This shows the surjection of $H^{n-3}((\Omega_{Y} \otimes (\omega_{Y}(D)))) \to H^{n-3}(\Omega_{D} \otimes \omega_D)$. \newline

Now let us finish the proof for $(2)$.
In the case (iii) of Theorem \ref{Fano}, $\omega_{Y_j}(D)$ is ample and $Y_j$ satisfies Bott-vanishing theorem, we have $H^{n-2}(\Omega_{Y_j} \otimes (\omega_{Y_j}(D))) = 0$. In the case (ii) and $d \geq 4$, we have that $D \in |\omega_Y^{-1}|$ and hence $h^{n-2}(\Omega_{Y_j} \otimes (\omega_{Y_j}(D))) = h^{n-2}(\Omega_{Y_j}) = h^{1,n-2}(Y_j) = 0$. This proves part $(2)$ of the lemma.  \newline
Now we prove part $(1)$
We show that $h^2(T_V) = 1$. We have an exact sequence
\begin{equation*}
    0 \to  H^{n-2}(\Omega_V/\tau_V \otimes \omega_V)  \to \bigoplus_{j = 1}^2 H^{n-2}(\Omega_{Y_j} \otimes (\omega_{Y_j}(D))) \to H^{n-2}(\Omega_{D} \otimes \omega_D) 
\end{equation*}
We show that the kernel is one-dimensional. In this case $D  = \omega_Y^{-1}$ and $n = 3$. So the above exact sequence simplifies to 
\begin{equation}\label{Tziolas exact sequence for Calabi-Yau}
    0 \to  H^{1}(\Omega_V/\tau_V \otimes \omega_V)  \to \bigoplus_{j = 1}^2 H^{1}(\Omega_{Y_j})  \to H^{1}(\Omega_{D}) 
\end{equation}
Since we have by assumption $h^{1,1}(Y_j) = 1$, if $Y_1 = Y$ is chosen to be general, then by result of Mukai, (see \cite{M88} or \cite{B02}, Theorem (page $1$) and Corollary $4.1$) the anticanonical $K3$ surface $D$ has $h^{1,1}(D) = 1$. Hence the above map cannot be injective and since it is not the zero map we have a one dimensional kernel. \QEDA

\section{Projective smoothing of union of Fano threefolds}\label{specializing to Fano threefolds}

\textcolor{black}{In this section we specialize to the case of Fano-threefolds. Due to the classification results of smooth Fano threefolds by Iskovskikh-Mori-Mukai, we are able to prove much more precise results and identify all the concrete cases of projective degeneration of smooth varieties into union of smooth Fano-threefolds that we obtain by our general theorems of the previous section. In each case we crucially use that the smoothing comes as the image of a $2:1$ morphism deforming to an embedding to find out the Picard group of the general fibre. To denote the deformation types, we use the notation according to \url{https://www.fanography.info/}.} Further we use the additive notation to denote line bundles and divisors in this section.

\smallskip

\textcolor{black}{In the following Corollary we give a list of cases where a union of anticanonically embedded Fano threefolds intersect along a del-Pezzo surface and are smoothable into smooth Fano threefolds which are themselves embedded by a sublinear series of the bi-anticanonical bundle. We also find out what the general smoothing fibres are.}

\begin{corollary}{\label{Fano-threefolds}}
 Let $Y \hookrightarrow \mathbb{P}^N$ be a anticanonically embedded smooth projective Fano threefold with index at least two. By Iskovskikh-Mori-Mukai's classification and notation according to \url{https://www.fanography.info/}, these are the families $1.12, 1.13, 1.14, 1.15, 1.16, 1.17$. Let $D \in  |-\displaystyle\frac{1}{2}K_Y|$ if index is $2$ and $D \in  |-\displaystyle\frac{1}{3}K_Y|$ or $D \in  |-\displaystyle\frac{2}{3}K_Y|$ if index is $3$, be a smooth divisor, which is a del-Pezzo surface. Then if $Y_1$ is the subvariety $Y$ of $\mathbb P^N$ and $Y_2$ is a subvariety of $\mathbb P^N$ that corresponds to a general element of $\mathcal{F}_{D,Y}$, then {the scheme theoretic intersection of $Y_1 and Y_2$ is $D$ and} {$V = Y_1 \bigcup_D Y_2$} is  smoothable inside $\mathbb{P}^N$ into \textbf{smooth varieties $X_t$ with ample anticanonical bundle or smooth Fano varieties}. The $Y_i, D, N$ and the general smooth fibres are listed below.   

\begin{itemize}
    \item[(1)]   If $Y_i$ belongs to type $1.12$, $D \in |-\displaystyle\frac{1}{2}K_{Y_i}|$, the smoothing fibres $X_t$ belongs to type $1.2$ (because $-K_{X_t}$ is base point free and $-K_{X_t}^3 = 4$) and are embedded in $\mathbb{P}^{10}$ by a sublinear series of $-2K_{X_t}$ of codimension $4$.
\vspace{0.2cm}
   
    \item[(2)]   If $Y_i$ belongs to type $1.13$, $D \in |-\displaystyle\frac{1}{2}K_{Y_i}|$, the smoothing fibres $X_t$ belongs to type $1.3$ (because $-K_{X_t}$ is very ample and in fact projectively normal and $-K_{X_t}^3 = 6$) and are embedded in $\mathbb{P}^{18}$ by a sublinear series of $-2K_{X_t}$ of codimension $6$.
    
    \item[(3)] If $Y_i$ belongs to type $1.14$, $D \in |-\displaystyle\frac{1}{2}K_{Y_i}|$, the smoothing fibres $X_t$ belongs to type $1.4$ (because $-K_{X_t}$ is very ample and in fact projectively normal and $-K_{X_t}^3 = 8$) and are embedded in $\mathbb{P}^{18}$ by a sublinear series of $-2K_{X_t}$ of codimension $6$.
    
    \item[(4)] If $Y_i$ belongs to type $1.15$, $D \in |-\displaystyle\frac{1}{2}K_{Y_i}|$, the smoothing fibres $X_t$ belongs to type $1.5(b)$ (by the explicit description of $1.5(b)$) with $-K_{X_t}^3 = 10$ and are embedded in $\mathbb{P}^{22}$ by a sublinear series of $-2K_{X_t}$ of codimension $7$
    
    \item[(5)] 
    \begin{itemize}
        \item[(a)]  If $Y_i$ belongs to type $1.16$, $D \in |-\displaystyle\frac{1}{3}K_{Y_i}|$, the smoothing fibres $X_t$ belongs to type $1.14$ (since $h^1(\Omega_{X_t}) = h^1(\Omega_{Y_i}) = 1$) with $-K_{X_t}^3 = 32$ and are embedded in $\mathbb{P}^{29}$ by a sublinear series of $-\displaystyle\frac{3}{2}K_{X_t}$ of codimension $14$.
       \item[(b)] If $Y_i$ belongs to type $1.16$, $D \in |-\displaystyle\frac{2}{3}K_{Y_i}|$, the smoothing fibres $X_t$ belongs to type $1.2 (b)$ (by the explicit description of $1.2(b)$) with $-K_{X_t}^3 = 4$ and are embedded in $\mathbb{P}^{22}$ by a sublinear series of $-3K_{X_t}$ of codimension $5$
    \end{itemize}

\item[(6)] If $Y_i$ belongs to type $2.32$, $D \in |-\displaystyle\frac{1}{2}K_{Y_i}|$, the smoothing fibres $X_t$ belongs to type $2.6(b)$ (by the explicit description of $2.6(b)$) with $-K_{X_t}^3 = 12$ and are embedded in $\mathbb{P}^{22}$ by a sublinear series of $-2K_{X_t}$ of codimension $8$

\item[(7)] If $Y_i$ belongs to type $2.35$, $D \in |-\displaystyle\frac{1}{2}K_{Y_i}|$, the smoothing fibres $X_t$ belongs to type $2.8$ (by the explicit description of $2.8$) with $-K_{X_t}^3 = 14$ and are embedded in $\mathbb{P}^{30}$ by a sublinear series of $-2K_{X_t}$ of codimension $9$.

\item[(8)] If $Y_i$ belongs to type $3.27$, $D \in |-\displaystyle\frac{1}{2}K_{Y_i}|$, the smoothing fibres $X_t$ belongs to type $3.1$ (by the explicit description of $3.1$) with $-K_{X_t}^3 = 12$ and are embedded in $\mathbb{P}^{26}$ by a sublinear series of $-2K_{X_t}$ of codimension $7$.    
\end{itemize}

In all cases the simple normal crossing Fano varieties lie in the boundary of Hilbert scheme of non-prime smooth Fano threefolds. In cases $(1)-(5)$ the Picard rank of the general smoothing is one but the hyperplane section does not generate the Picard group. If $P$ denote the generator of the Picard group and $A$, the hyperplane section then in cases $(1)-(4)$, $A = 2P$, in case $5(a)$, $A = -\displaystyle\frac{3}{2}P$ and in case $5(b)$, $A = 3P$.  In cases $(6)-(7)$ the Picard rank of the general smoothing is two while in case $(8)$, the Picard rank of the general smoothing is $3$.

\end{corollary}

\noindent\textit{Proof.} The proof follows from Theorem \ref{Fano}, $(c)$(i). Let $Y \hookrightarrow \mathbb{P}^N$ be a anticanonically embedded Fano threefold with index $2$ or $3$. By Iskovskikh-Mori-Mukai's classification and notation according to \url{https://www.fanography.info/}, these are $1.12, 1.13, 1.14, 1.15, 1.16, 1.17, 2.32, 2.35, 3.27$. If the index is $2$, we can take $L$ such that $2L=K_Y$, $H = -K_Y$ while if index is $3$, we can take $L$ such that $3L = K_Y$, $H = -K_Y$ or $L$ such that \textcolor{black}{$3L=2K_Y$}, $H = -K_Y$ in Theorem \ref{Fano}, part $(c)$(i).  The line bundle $-\displaystyle\frac{1}{2}K_Y$ or $-\displaystyle\frac{1}{3}K_Y$ or $-\displaystyle\frac{2}{3}K_Y$ is ample and base point free and contains a smooth divisor $D$.  Note that in all the cases listed, either $H + L = -\displaystyle\frac{1}{2}K_Y$ or $H + L = -\displaystyle\frac{2}{3}K_Y$ or $H + L = -\displaystyle\frac{1}{3}K_Y$ which are base point free. Moreover in each of the cases, $h^0(-K_Y) \geq 8$. These will give us the existence of union of anticanonically embedded Fano-threefolds intersecting along del-Pezzo surfaces which are smoothable into smooth Fano threefolds. To see what the general smoothing fibers $X_t$ are, we consider the Fano double cover $\pi: X \to Y$ branched along a smooth curve $B$ in the linear system $|2D| = |-2L|$. The smoothing fibers $X_t$ are deformations of $X$ along some one parameter family. Then in case $(1)$, we note that $-K_X$ is base point free and $-K_{X}^3 = 4$ and in cases $(2)$ and $(3)$ we note that $-K_X$ is projectively normal and hence in particular very ample and $-K_{X}^3 = 6$ and $-K_{X}^3 = 8$ respectively. Hence the same is true for $-K_{X_t}$. Now we can conclude the deformation type of $X_t$ since there is only one deformation type of Fano-threefolds with the above property of $-K_{X_t}$ and the given value of $-K_{X_t}^3$. In case $(5)(a)$, we first note that $-K_{X_t}^3 = 32$. Then we calculate $h^{1,1}(X) = h^1(\Omega_X)$ as follows: $H^1(\Omega_X^1) = H^1(\pi_*\Omega_X^1) = H^1(\Omega_Y) \bigoplus H^1(\Omega_Y(\log (B)) \otimes L)$. Now using the exact sequence 
$$ 0 \to \Omega_Y \otimes L \to \Omega_Y(\log (B)) \otimes L \to L|_B \to 0 $$ we have 
$$H^1(\Omega_Y \otimes L) \to H^1(\Omega_Y(\log (B)) \otimes L) \to H^1(L|_B)$$
The leftmost term is zero due to Nakano vanishing theorem. We have 
$$0 \to 3L \to L \to L|_B \to 0$$ 
Since $-L$ is ample and $Y$ is Fano, we have that $H^1(L) = H^2(3L) = 0$ and hence we see that $H^1(L|_B) = 0$. This implies $H^1(\Omega_Y(\log (B)) \otimes L) = 0$ and consequently $h^1(\Omega_X) = h^1(\Omega_Y) = 1$ since $Y$ is of type $1.16$. This implies $h^1(\Omega_{X_t}) = 1$ since Hodge numbers are constant along a smooth family. Now once again we conclude by noting there is only deformation type of Fano threefold with $-K_{X_t}^3 = 32$ and $h^1(\Omega_{X_t}) = 1$. In the rest of the cases, once again the smoothing fibers $X_t$ are deformations of the Fano double cover $\pi: X \to Y$, but in all these cases, a general deformation of a such a double cover is once again a double cover which we can conclude from the list in \url{https://www.fanography.info/}. We now show that the smoothing fibers are embedded by a sublinear series of $|-mK_{X_t}|$ with $m$ and codimension as stated in the Theorem. Note that $K_X = \pi^*(K_Y-L)$ and 
consider the composition  
\[
\begin{tikzcd}
X \arrow[dr, "\varphi"] \arrow[d, "\pi"] & \\
Y \arrow[r, hook, "|-K_Y|"] & \mathbb{P}^N
\end{tikzcd}
\]
$X_t \xhookrightarrow{\varphi_t} \mathbb{P}^N$ is obtained as the image of an embedding $\varphi_t$ which is a deformation of $\varphi$ along a one parameter family $T$. Let $A = \varphi^*(\mathcal{O}_{\mathbb{P}^N}(1)) = \pi^*(-K_Y)$. Plugging in the value of $L$ in each case in $K_X = \pi^*(K_Y-L)$, we see that $A = -2K_X$ in cases $(1)-(4)$ and $(6)-(8)$ while $A = -\displaystyle\frac{3}{2}K_X$ in case $(5)(a)$ and $A = -3K_X$ in case $(5)(b)$. The one parameter family of smooth fibres $X_t$ that smoothes the normal crossing $Y_1 \bigcup_D Y_2$ are embedded inside $\mathbb{P}^N$ by some sublinear series of a line bundle $A_t$ which is a deformation of $A$. But $X$ is regular since $Y$ is regular and $H^1(L) = 0$. This implies $A_t = -mK_{X_t}$ with $m$ as mentioned in the Theorem. Now note that $h^0(A) = h^0(\pi^*(-K_Y)) = h^0(-K_Y)+h^0(-K_Y+L) = N+1+h^0(-K_Y+L)$. Also for $i = 1,2,3$, $h^i(A) = 0$. This implies that $h^0(A_t) = h^0(A) = N+1+h^0(-K_Y+L)$. Hence $X_t$ is embedded by a sublinear series of $-mK_{X_t}$ of codimension $h^0(-K_Y+L)$. Now that we have identified the hyperplane section of the smoothing fibre and we also know the generators of its Picard groups from the classification in \url{https://www.fanography.info/}, the rest of the statement follows.  \QEDA \newline

\textcolor{black}{In the following Corollary we give a list of cases where a union of anticanonically or bi-anticanonically embedded Fano threefolds intersect along a $K3$ surface and are smoothable into smooth Calabi-Yau threefolds.} 

\begin{corollary}{\label{Calabi-Yau threefolds}}
Let $Y \hookrightarrow \mathbb{P}^N$ be an embedding of a smooth projective Fano threefold with $h^{1,1}(Y) = 1$. Let $D \in  |-K_Y|$ be a smooth divisor, which is a $K3$ surface.
\begin{itemize}
    \item[(i)]   If $h^0(-K_Y) \geq 8$, then assume that the embedding $Y \hookrightarrow \mathbb{P}^N$ is induced by $|-K_Y|$. These are the families $1.5-1.10$ and $1.12-1.17$.  
    \item[(ii)]  If $h^0(-K_Y) \leq 7$, then assume that the embedding $Y \hookrightarrow \mathbb{P}^N$ is induced by $|-2K_Y|$. These are the families $1.1-1.4$ and $1.11$. 
\end{itemize}
Then 
    \begin{itemize}
    \item[(1)] 
     \begin{itemize}   
        \item[(a)] If $Y_1$ is the subvariety $Y$ of $\mathbb P^N$ and $Y_2$ is a subvariety of $\mathbb P^N$ that corresponds to a general element of $\mathcal{F}_{D,Y}$, then {the scheme theoretic intersection of $Y_1 and Y_2$ is $D$ and} {$V = Y_1 \bigcup_D Y_2$} is  smoothable inside $\mathbb{P}^N$ into \textbf{smooth Calabi-Yau threefolds} embedded by a sublinear series of codimension one of a very ample line bundle in case (i) and of codimension $h^0(-K_Y)$ in case (ii). Moreover, the Hilbert point of $V$ inside $\mathbb{P}^N$ is smooth. 
         
        \item[(b)] For families $1.5-1.10$, the fibre $\mathcal{F}_{D}$ is irreducible, i.e, $\mathcal{F}_{D} = \mathcal{F}_{D,Y}$  and hence all normal crossing varieties $V = Y_1 \bigcup_D Y_2$, embedded inside $\mathbb{P}^N$ where $Y_2$ is a deformation of $Y_1 = Y$, that intersect along the subscheme $D$ form an irreducible family and a general such $V$ is smoothable inside $\mathbb{P}^N$. 
    \end{itemize}
    
    \item[(2)] The Picard group of a general Calabi-Yau threefold of the component of the Hilbert scheme found in part $(1)$ is $\mathbb{Z}$. In the case $(i)$, if the index of $Y$ is one, i.e, if $Y$ belongs to families $1.5-1.10$, the general embedded smoothing is a prime Calabi Yau threefold while if the index of $Y$ is atleast two, i.e, if $Y$ belongs to families $1.12-1.17$, the general embedded smoothing is a non-prime Calabi Yau threefold and the hyperplane section is a multiple (which equals the index of the Fano) of the generator of the Picard group. In case $(ii)$ (and hence the index of $Y$ is one), the general embedded smoothing is a non-prime Calabi Yau threefold and the hyperplane section is twice of the generator of the Picard group

\end{itemize}

\end{corollary}

\noindent\textit{Proof.} Part $(1)(a)$ follows from Theorem ~\ref{Fano}, $(c)$(ii). If $Y \hookrightarrow \mathbb{P}^N$ be an embedded Fano threefold with Picard rank $\rho = 1$ and as listed in the Theorem, then we can choose $L = K_Y$, $D \in  |-K_Y|$ smooth and $H = -K_Y$ in case (i) and $L = K_Y$, $D \in  |-K_Y|$ smooth and $H = -2K_Y$ in case (ii) and apply Theorem ~\ref{Fano}, $(c)$ (ii) to get the existence and smoothing of the union of the embedded Fano-threefolds into Calabi-Yau threefolds. We now prove the statement on the codimension of the sublinear series with which the general fibers are embedded. Consider a Calabi-Yau double cover of $\pi:X \to Y$ branched along $-2L = -2K_Y$. Let $A = \pi^*(\mathcal{O}_{Y}(1))$ and consider the composition  
\[
\begin{tikzcd}
X \arrow[dr, "\varphi"] \arrow[d, "\pi"] & \\
Y \arrow[r, hook] & \mathbb{P}^N
\end{tikzcd}
\]
$X_t \xhookrightarrow{\varphi_t} \mathbb{P}^N$ is obtained as the image of an embedding $\varphi_t$ which is a deformation of $\varphi$ along a one parameter family $T$. Now note that $\varphi^* (\mathcal{O}_{\mathbb{P}^N}(1)) = \pi^*\mathcal{O}_Y(1)$. Then in case (i), $A = \pi^*(-K_Y)$ and in case (ii), $A = \pi^*(-2K_Y)$. Note that the one parameter family of smooth fibers $X_t$ that smoothes the normal crossing $Y_1 \bigcup_D Y_2$ are embedded inside $\mathbb{P}^N$ by some sublinear series of a line bundle $A_t$ which is a deformation of $A$. Since $X_t$ is Calabi-Yau, we have that for $i = 1,2,3$, $H^i(A_t) = 0$ and hence $h^0(A_t) = h^0(A)$. In case (i), $h^0(A_t) = h^0(A) = h^0(\pi^*(-K_Y)) = h^0(-K_Y)+h^0(\mathcal{O}_Y) = N+1+1$. So in this case, the codimension of the sublinear series that gives the embedding is one. In case (ii), $h^0(A_t) = h^0(A) = h^0(\pi^*(-2K_Y)) = h^0(-2K_Y)+h^0(K_Y) = N+1+h^0(K_Y)$. So in this case, the codimension of the sublinear series that gives the embedding is $h^0(K_Y)$. 
$(1)(b)$ follows from the penultimate statement of Theorem $7$, \cite{CLM93} and Theorem \ref{intrinsic}, $(2)$.  
\newline 
\noindent\textit{Proof (2).} Consider the Calabi-Yau double cover $\pi: X \to Y$ branched along $B \in |-2K_Y|$. By \cite{Par91}, Lemma $4.2$, we have $H^1(\Omega_X^1) = H^1(\pi_*\Omega_X^1) = H^1(\Omega_Y) \bigoplus H^1(\Omega_Y(\log (B)) \otimes K_Y)$. Now using the exact sequence 
$$ 0 \to \Omega_Y \otimes K_Y \to \Omega_Y(\log (B)) \otimes K_Y \to K_Y|_B \to 0 $$ we have 
$$H^1(\Omega_Y \otimes K_Y) \to H^1(\Omega_Y(\log (B)) \otimes K_Y) \to H^1(K_Y|_B)$$
The leftmost term is zero due to Nakano vanishing theorem. We have 
$$0 \to 3K_Y \to K_Y \to K_Y|_B \to 0$$ 
and hence we see that $H^1(K_Y|_B) = 0$. This implies $H^1(\Omega_Y(\log (B)) \otimes K_Y) = 0$. Now since $H^1(\mathcal{O}_Y) = H^1(\mathcal{O}_X) = 0$, the pullback map between $\pi^*: \textrm{Pic}(Y)  \to \textrm{Pic}(X) $ is induced by the map 
$$H^1(\Omega_Y) \to H^1(\Omega_Y) \bigoplus H^1(\Omega_Y(\log (B)) \otimes K_Y)$$
Therefore $H^1(\Omega_Y(\log (B)) \otimes K_Y) = 0$ implies $\pi^*$ is an isomorphism and $h^{1,1}(X) = 1$. Since the general smoothing $X_t$ of $V$ is a deformation of $X$, we have that $h^{1,1}(X_t) = 1$. On the other hand, by \cite[Corollary 4.18]{Debarre} and \cite[Theorem 0.1]{KMM},  $Y$ is simply connected, so, by the Lefschetz theorem  for homotopy groups (see \cite[Theorem 3.1.21]{positivity}), so is $B$, where $B$ is the ample branch divisor of the double cover $X$. Now if $R \subset X$ denotes the ramification divisor, $R \cong B$ and hence $R$ is simply connected. Since $R$ is ample, once again by the Lefschetz theorem  for homotopy groups, so is $X$. Therefore $X$ is simply connected and hence so is $X_t$. Therefore, $\textrm{Pic}(X_t) = \mathbb{Z}$. We  prove the rest of the statement for case $(i)$. The proof for case $(ii)$ is similar to proof for families $1.12-1.17$. Now suppose that index of $Y$ is one or $-K_Y$ generates $\textrm{Pic}(Y)$. We want to show that the hyperplane section generates the Picard group of $X_t$.   
Consider the composition  
\[
\begin{tikzcd}
X \arrow[dr, "\varphi"] \arrow[d, "\pi"] & \\
Y \arrow[r, hook, "|-K_Y|"] & \mathbb{P}^N
\end{tikzcd}
\]
$X_t \xhookrightarrow{\varphi_t} \mathbb{P}^N$ is obtained 
as the image of an embedding $\varphi_t$ which is a deformation of $\varphi$ along a one parameter family $T$. Now note that $\varphi^* (\mathcal{O}_{\mathbb{P}^N}(1)) = \pi^*(-K_Y)$ and since $\pi^*$ is an isomorphism and $-K_Y$ generates $\textrm{Pic}(Y)$, we have that $\pi^*(-K_Y)$ maps to the generator of $H^1(\Omega_{X})$. Hence we have that the map 
$$ H^1(\Omega_{\mathbb{P}^N}) \to H^1(\Omega_{X}) $$ is an isomorphism. Since this map globalizes along the one parameter family $T$, and both $H^1(\Omega_{\mathbb{P}^N})$ and $H^1(\Omega_{X_t})$ remains constant for $t \in T$, for a general $t \in T$,
$$ H^1(\Omega_{\mathbb{P}^N}) \to H^1(\Omega_{X_t}) $$ is an isomorphism. Hence $X_t$ is a prime Calabi-Yau threefold. Now if the index $m$ of $Y$ is at least $2$, then $-K_Y = mW$ where $W$ generates the Picard group of $Y$.  Consider the Calabi-Yau double cover $\pi: X \to Y$ and the line bundle $A = \pi^*(-K_Y)$.  Then $A = mA'$ where $A' = \pi^*(W)$. The hyperplane section $A_t$ of a general smoothing is a deformation of $A$.  Now note that $H^2(\mathcal{O}_X) = H^2(\mathcal{O}_Y) \bigoplus H^2(K_Y) = 0$. Hence $A'$ deforms to a line bundle $A'_t$ on $X_t$. Since $X$ is regular, this implies that $A_t = mA'_t$. \QEDA

\begin{remark}\label{Fano and Calabi Yau}

Using Corollary ~\ref{Fano-threefolds} and Corollary ~\ref{Calabi-Yau threefolds} we exhibit two different types of unions of anticanonically embedded smooth Fano varieties with Picard rank $1$ and index $\rho \geq 2$. There are five such families $1.12-1.16$. So for each such family we construct a union of anticanonically embedded Fano varieties intersecting along a del-Pezzo surface (respectively a $K3$ surface) which are smoothable into a smooth Fano (respectively a Calabi-Yau threefold) inside the projective space $\mathbb{P}^{h^0(-K_Y)}$. For $1.16$ we have in fact $3$ different types of unions 

\end{remark}

\begin{remark}\label{deviation from Kulikov}
  We point out an interesting difference in the higher dimensional case. By the result of Kulikov (see \cite{Ku77}), in a type II degeneration of $K3$ surfaces with two irreducible components, both components must be rational, while in Theorem \ref{Calabi-Yau threefolds}, for families $1.5, 1.7, 1.12, 1.13$ (respectively for familes $1.1, 1.2, 1.3, 1.4, 1.11$) we get examples of degenerations of polarized Calabi-Yau varieties into snc union of anticanonically (respectively bi-anticanonically) embedded Fano varieties, neither of whose components are rational.  
\end{remark}

In the following Corollary we give a list of cases where a union of bi-anticanonically embedded Fano threefolds intersect along canonically embedded surfaces and are smoothable into smooth threefolds with ample canonical bundle. Note that by Theorem \ref{Fano}, we require the Fano threefolds to satisfy Bott vanishing Theorem. In \cite{Tot23}, the author has classified all smooth Fano threefolds satisfying Bott vanishing theorem.

\begin{corollary}\label{Fano with Bott}
 Let $Y \hookrightarrow \mathbb{P}^N$ be a bi-anticanonically embedded smooth projective Fano threefold satisfying Bott vanishing theorem. According to \cite{Tot23}, and by Iskovskikh-Mori-Mukai's classification and notation according to \url{https://www.fanography.info/}, these are the families $1.17$, $2.33-2.36$, $3.25-3.31$, $4.9-4.12$, $5.2-5.3$ (toric), $2.26$, $2.30$, $3.15-3.16$, $3.18-3.24$, $4.3-4.8$, $5.1$, $6.1$ (non-toric). Let $D \in  |-2K_Y|$ be a smooth divisor, which is a surface of general type. Then
\begin{itemize}
    \item[(1)]  If $Y_1$ is the subvariety $Y$ of $\mathbb P^N$ and $Y_2$ is a subvariety of $\mathbb P^N$ that corresponds to a general element of $\mathcal{F}_{D,Y}$, then {the scheme theoretic intersection of $Y_1 and Y_2$ is $D$ and} {$V = Y_1 \bigcup_D Y_2$} is  smoothable inside $\mathbb{P}^N$ into \textbf{smooth threefolds of general type } embedded by a sublinear series of codimension $1$ of their bicanonical linear system. Moreover, the Hilbert point of $V$ inside $\mathbb{P}^N$ is smooth.
    
    \item[(2)] The Picard group of a general threefold of general type  of the component of the Hilbert scheme found in part $(1)$ is isomorphic to the Picard group of $Y$. If $\textrm{Pic}(Y) = \mathbb{Z}$ then the hyperplane section is a multiple (which equals the index of $Y$) of the generator of the Picard group.  

\end{itemize}

\end{corollary}

\noindent\textit{Proof $(1)$.} Part $(1)$ follows from Theorem \ref{Fano}, part $(c)$(iii), by taking $L = 2K_Y$, $D \in |-2K_Y|$ smooth and $H = -2K_Y$. We show that the general smoothing fibres of $V$ are embedded by a sublinear series of codimension $1$ of their bicanonical linear series. Consider a double cover of $\pi:X \to Y$ branched along $-2L = -4K_Y$. Then $K_X = \pi^*(-K_Y)$ and is hence ample and base point free. Let $A = \pi^*(-2K_Y) = 2K_X$. Note that the one parameter family of smooth fibres $X_t$ that smoothes the normal crossing $Y_1 \bigcup_D Y_2$ are embedded inside $\mathbb{P}^N$ by some sublinear series of a line bundle $A_t$ which is a deformation of $A = 2K_X$. But $X$ is regular since $Y$ is regular and $H^1(2K_Y) = 0$. This implies $A_t = 2K_{X_t}$. Now note that $h^0(A) = h^0(\pi^*(-2K_Y)) = h^0(-2K_Y)+h^0(\mathcal{O}_Y) = h^0(-2K_Y)+1$. Also for $i = 1,2,3$, $h^i(\pi^*(-2K_Y)) = 0$. This implies that $h^0(A_t) = h^0(A) = (N+1)+1$. Hence $X_t$ is embedded by a sublinear series of $2K_{X_t}$ of codimension $1$. \\
\noindent\textit{Proof $(2)$.} The proof is very similar to part $(2)$ of Corollary \ref{Calabi-Yau threefolds}. Once again we have an exact sequence
$$H^1(\Omega_Y \otimes 2K_Y) \to H^1(\Omega_Y(\log (B)) \otimes 2K_Y) \to H^1(2K_Y|_B)$$
where $B \in |-4K_Y|$. By Nakano vanishing theorem we have the vanishing of $H^1(\Omega_Y \otimes 2K_Y)$ and the vanishing of the right hand side follows from the defining exact sequence of $B$. This shows that 
$$H^1(\Omega_Y) \cong H^1(\Omega_X)$$ 
is an isomorphism. 
But now we know from the classification that the $\textrm{Pic}(Y) \otimes \mathbb{C} = H^1(\Omega_Y)$. Combining with the facts that $\pi^*$ is an injective map, both $H^1(\mathcal{O}_Y) = H^1(\mathcal{O}_X) = 0$, we have that 
$$\textrm{Pic}(Y) \cong \textrm{Pic}(X)$$
under $\pi^*$. 
Now since $h^1(\Omega_Y) \cong h^1(\Omega_X)$ we have that the general smoothing $X_t$ has Picard number less than or equal to the Picard number of $Y$. Also, arguing similarly as in the proof of Corollary \ref{Calabi-Yau threefolds}, $(2)$, we find that $\textrm{Pic}(X_t)$ is torsion free. Now we show that $X_t$ has Picard number exactly equal to the Picard number of $Y$. Let the generators of $\textrm{Pic}(Y)$ be denoted by $L_i$ and the generators of $\textrm{Pic}(X)$ be denoted by $\widetilde{L}_i = \pi^* L_i$. Since $L_i$ generates $\textrm{Pic}(Y)$ freely, we have that $\widetilde{L}_i$ generates $\textrm{Pic}(X)$ freely. Now since $H^2(\mathcal{O}_X) = 0$ all the $\widetilde{L}_i$'s deform to line bundles $L_{it}$ on the general fibre $X_t$. Now suppose that $L_{it} = \Sigma_{j\neq i}a_jL_{jt}$ for some $a_j$. Then consider the line bundle $\widetilde{L}_{i} - \Sigma_{j\neq i}a_j\widetilde{L}_{j}$ on $X$. This is the limit of $L_{it} - \Sigma_{j\neq i}a_jL_{jt}$ as $t \to 0$. Since numerical chern classes are constant along a smooth family, $c_1(\widetilde{L}_{i} - \Sigma_{j\neq i}a_j\widetilde{L}_{j}) = 0$. But now since $H^1(\mathcal{O}_{X}) = 0$, we have $\widetilde{L}_{i} - \Sigma_{j\neq i}a_j\widetilde{L}_{j} = 0$ which contradicts the fact that $\widetilde{L}_i$'s freely generate $\textrm{Pic}(X)$.  
The proof of the assertion when $\textrm{Pic}(Y) = \mathbb{Z}$ is exactly same as the corresponding proof in Corollary \ref{Calabi-Yau threefolds}. \QEDA

\section{Construction of smoothable polarized Calabi-Yau and stable pairs}

As we already mentioned, one of the consequences of an embedded smoothing of $V = Y_1 \bigcup_D Y_2 \subset \mathbb{P}^N$ is that it automatically gives us a smoothing of the pair $(V, \Delta = cH)$ where $H$ is any hyperplane section of the embedding. In particular choosing $H$ to be general and $c$ to be a suitable rational value, we obtain smoothable slc polarized Calabi-Yau pairs from Corollary \ref{Fano-threefolds} and stable pairs  from Corollaries \ref{Fano-threefolds}, \ref{Calabi-Yau threefolds} and \ref{Fano with Bott}. We state the precise results in this section and list some of the examples explicitly in Examples \ref{examples of log pairs from Fano} and \ref{examples of log pairs from CY and gen type}. First we recall some definitions from \cite{KX20}: 

\begin{definition}\label{stable and polarized Calabi-Yau pairs}
\begin{itemize}
    \item[(1)] A pair $(X, \Delta)$ consists of a reduced equidimensional variety $X$ and an effective $\mathbb{Q}-$ divisor $\Delta$ none of whose irreducible components are contained in singular locus of $X$.
    \item[(2)] A pair $(X, \Delta)$ is stable if it is semi-log canonical, $X$ is projective and $K_X+\Delta$ is an ample $\mathbb{Q}-$ divisor. A morphism $f: (X, \Delta) \to S$ with $S$ normal, is called stable if every fibre over closed points of $S$ is stable, $K_{X/S}+\Delta$ is $\mathbb{Q}-$ Cartier and $f-$ ample.
    \item[(3)]  A pair $(X, \Delta)$ is a Calabi-Yau pair, if it is semi-log canonical, $X$ is proper and $K_X+\Delta$ is $\mathbb{Q}-$ linearly equivalent to $0$. We write $\Delta = cD$ where $D$ is a $\mathbb{Z}-$ divisor and consider $c$ to be fixed in our moduli problem. 
    \item[(4)] A polarized Calabi-Yau pair $(X, \Delta)$ is a Calabi-Yau pair, plus an ample $\mathbb{Q}-$ Cartier divisor $H$ such that $(X, 
    \Delta+\epsilon H)$ is semi-log-canonical and consequently a stable pair for $0 < \epsilon << 1$. The latter holds iff $H$ does not contain any of the slc centers of $(X, \Delta)$. 
    \item[(5)] A stable family of polarized Calabi-Yau pairs over a normal base scheme $S$ consists of a flat proper morphism $f: X \to S$, a $\mathbb{Q}-$ divisor $\Delta$ and a $\mathbb{Q}-$ Cartier divisor $H$ on $X$ such that $K_{X/S}+\Delta$ is $\mathbb{Q}-$ Cartier and all fibers $(X_s, \Delta_s, H_s)$ are polarized Calabi-Yau pairs.  
\end{itemize}    
\end{definition}

It is known that in characteristic $0$, there is a coarse moduli space of stable pairs $\textrm{\textbf{SP}}$ which is separated and satisfies the valuative criterion of properness and whose connected components are projective. It is known that polarized Calabi-Yau pairs have a coarse moduli space $\textrm{\textbf{PCY}}$ whose irreducible components are projective.

\begin{theorem}\label{pairs main}
In the situation of Theorem \ref{intrinsic}, further assume that $Y$ is regular. If $d, m$ and $n$ are integers with $d \geq n$ such that 
 $$d(K_Y+D)+mnH = 0 (\textrm{respectively} \hspace{0.1cm} d(K_Y+D)+mnH \textrm{is ample})$$ where $H$ is the hyperplane section of the embedding $i: Y \hookrightarrow \mathbb{P}^N$, then for each smoothable $V \subset \mathbb{P}^N$ constructed in Theorem \ref{intrinsic} part $(1)$, and a general hyperplane section $B_m \in \mathcal{O}_V(m)$ there exists a Gorenstein one parameter family $(\mathcal{V}, \displaystyle\frac{n}{d}\mathcal{B}_m) \to T$ of pairs such that
    \begin{itemize}
        \item[(a)] Both $K_{\mathcal{V}/T}$ and $\mathcal{B}_m$ are Cartier.
        \item[(b)] $(\mathcal{V}_0, \displaystyle\frac{n}{d}\mathcal{B}_{m0}) = (V, \displaystyle\frac{n}{d}B_m)$ is a Calabi-Yau pair (respectively a stable pair) where $B_m = H_{m1} \bigcup_C H_{m2} \in |\mathcal{O}_V(m)|$ is a simple normal crossing Cartier divisor in $V$, such that $H_{mi} \in |\mathcal{O}_{Y_i}(m)|$ are smooth divisors in $Y_i$ and $C = H_{m1} \bigcap D = H_{m2} \bigcap D$ is smooth.
        \item[(c)] For $t \neq 0$,  $(\mathcal{V}_t, \displaystyle\frac{n}{d}\mathcal{B}_t) = (V_t, \displaystyle\frac{n}{d}B_{mt})$ is a Calabi-Yau pair (respectively a stable pair) where $V_t$ is smooth and $B_{mt} \subset V_t$ is a smooth Cartier divisor.
    \end{itemize}

\end{theorem}

\noindent\textit{Proof.} Let $\mathcal{V} \subset \mathbb{P}^N_T \to T$ be an embedded one parameter family over $T$, such that $\mathcal{V}_0 = V$ and $\mathcal{V}_t$ is smooth. This in particular implies that $\mathcal{O}_{\mathcal{V}}(1)$ is a relatively very ample line bundle and hence so is $\mathcal{O}_{\mathcal{V}}(m)$ for any $m$. Now consider a general hyperplane section $\mathcal{B}_m \in |\mathcal{O}_{\mathcal{V}}(m)|$ obtained from an embedding given by $\mathcal{O}_{\mathcal{V}}(m)$. We have that $\mathcal{B}_{m0} = B_m \in |\mathcal{O}_{V}(m)|$ is a simple normal crossing divisor $B = H_{m1} \bigcup_C H_{m2}$ such that $H_{mi} \in |\mathcal{O}_{Y_i}(m)|$ are smooth divisors in $Y_i$ and $C = H_{m1} \bigcap D = H_{m2} \bigcap D$ is smooth. Also for $t \neq 0$, $\mathcal{B}_{mt} = B_{mt} \in |\mathcal{O}_{V_t}(m)|$ is smooth. Now consider the pair $(\mathcal{V}, \displaystyle\frac{n}{d}\mathcal{B}_m)$. Since $\mathcal{B}_m$ is a Cartier divisor, $\displaystyle\frac{n}{d}\mathcal{B}_m$ is $\mathbb{Q}-$ Cartier. Since $V$ is Gorenstein, $V_t$ is Gorenstein for a general $t$, 
and hence after possibly shrinking $T$, the map $\mathcal{V} \to T$ is Gorenstein and hence $K_{\mathcal{V}/T}$ is Cartier (see \href{https://stacks.math.columbia.edu/tag/0C14}{Tag 0C14}). Now consider the central fibre $(V, \displaystyle\frac{n}{d}B_m)$. The normalization $\nu: V^{\nu} \to V$ is the disjoint union of $Y_1$ and $Y_2$ and the pair $(V^{\nu}, \displaystyle\frac{n}{d}\nu^{-1}B_m+D^{\nu})$ where $D^{\nu}$ is the inverse image of the double locus $D$ of $V$ is actually is the disjoint union of pairs $(Y_i, \displaystyle\frac{n}{d}H_{mi}+D)$ for $i = 1,2$, both of which are log canonical since $d \geq n$ and $H_{mi}$ and $D$ are smooth and do not contain any common components. Also the pairs $(V_t, \displaystyle\frac{n}{d}B_{mt})$ are lc since $B_{mt}$ is smooth and $d \geq n$. So all that is left to show is that $K_{V}+\displaystyle\frac{n}{d}B_m$ and $K_{V_t}+\displaystyle\frac{n}{d}B_{mt}$ are $\mathbb{Q}-$ linearly equivalent to $0$ ( resp. ample $\mathbb{Q}-$ Cartier divisors). Note that since $\mathbb{Q}-$ linear equivalence to the trivial divisor (resp.  ampleness) is an open condition, it is enough to prove the statement for $K_{V}+\displaystyle\frac{n}{d}B_m$. Since $V = Y_1 \bigcup_D Y_2$, we need to show that $K_{V}+\displaystyle\frac{n}{d}B_m|_{Y_i}$ is trivial (resp. ample) for each $i$. Note that $K_{V}|_{Y_i} = K_{Y_i}+D$ and hence $K_{V}+\displaystyle\frac{n}{d}B_m|_{Y_i} = K_{Y_i}+D+\displaystyle\frac{n}{d}H_{mi}$. Since $Y_1 = Y$, by assumptions, this is $\mathbb{Q}-$ Cartier and trivial (resp. ample ). Now since $Y_2$ is a general deformation of $Y_1$ fixing $D$ and $\mathbb{Q}-$ linear equivalence to the trivial divisor (resp. ampleness) is an open condition, the same conclusion holds for $Y_2$. \QEDA
\newline

\textcolor{black}{In the following Corollary we give a list of polarized Calabi-Yau pairs $(V,\Delta,\epsilon H)$ and stable pairs $(V,\Delta)$ where $V = Y_1 \bigcup_D Y_2$ is a union of two Fano threefolds intersecting along a del-Pezzo surface $D$, $H$ is anticanonical when restricted to each $Y_i$ and the pair is smoothable into polarized Calabi-Yau pairs $(V_t, \Delta_t, \epsilon H_t)$ or stable pairs $(V_t, \Delta_t)$ where $V_t$ is smooth.}

\begin{corollary}\label{log CY and gen type pairs from Fano}
The following smoothable Calabi-Yau pairs (resp. stable pairs) are obtained from Corollary \ref{Fano-threefolds} using Theorem \ref{pairs main} :
Let $d, m$ and $n$, $d \geq n$, be integers that satisfy $d= 2mn$ (resp. $d < 2mn$) in cases $(1)-(4), (6)-(8)$, $2d = 3mn$ (resp. $2d < 3mn$) in case $(5)(a)$ and $d = 3mn$ (resp. $d < 3mn$) in case $(5)(b)$ of Corollary \ref{Fano-threefolds}. Then for each smoothable $V = Y_1 \bigcup_D Y_2 \subset \mathbb{P}^N$, (where $Y_i$ are anticanonically embedded Fano threefolds both of the same deformation type corresponding to each case and $D$ is a subanticanonical del Pezzo surface) and a general  hyperplane section $B_m \in \mathcal{O}_V(m)$, there exists a Gorenstein one parameter family $(\mathcal{V}, \displaystyle\frac{n}{d}\mathcal{B}_m) \to T$ of pairs such that
    \begin{itemize}
        \item[(a)] Both $K_{\mathcal{V}/T}$ and $\mathcal{B}_m$ are Cartier.
        \item[(b)] $(\mathcal{V}_0, \displaystyle\frac{n}{d}\mathcal{B}_{m0}) = (V, \displaystyle\frac{n}{d}B_m)$ is a Calabi-Yau pair (resp. a stable pair) where $B_m = H_{m1} \bigcup_C H_{m2}$ is a simple normal crossing Cartier divisor in $V$, such that $H_{mi} \in |\mathcal{O}_{Y_i}(m)|$ are smooth divisors in $Y_i$ and $C = H_{m1} \bigcap D = H_{m2} \bigcap D$ is smooth.
        \item[(c)] For $t \neq 0$,  $(\mathcal{V}_t, \displaystyle\frac{n}{d}\mathcal{B}_t) = (V_t, \displaystyle\frac{n}{d}B_{mt})$ is a lc Calabi-Yau pair (resp. stable lc pair) where $V_t$ is a smooth Fano threefold and $B_{mt} \subset V_t$ is a smooth Cartier divisor.
   \end{itemize}

\end{corollary}

\begin{example}\label{examples of log pairs from Fano}
By taking another general hyperplane section $\mathcal{C}_m \in |\mathcal{O}_{\mathcal{V}}(m)|$, we get smoothable polarized Calabi-Yau pairs from Corollary \ref{log CY and gen type pairs from Fano}.  We describe explicitly some of the smoothable polarized CY obtained from Corollary \ref{log CY and gen type pairs from Fano}. 
\begin{itemize}
    \item[(1)] For $m = 1, n = 1, d = 2$ in cases $(1)-(4), (6)-(8)$, 
    \begin{itemize}
        \item[(a)]  in the central fibre we have a polarized Calabi-Yau pair \\ $(Y_1 \bigcup_D Y_2, \displaystyle\frac{1}{2}(H_{11} \bigcup_C H_{12})+\epsilon (H_{11}' \bigcup_{C'} H_{12}'))$ for $\epsilon < 1$, where $Y_i$ are Fano threefolds with index two, both of same deformation type i.e, of both are of type either $1.12$ or $1.13$ or $1.14$ or $1.15$ or $2.32$ or $2.35$ or $3.27$ respectively, $D \in |-\displaystyle\frac{1}{2}K_{Y_i}|$ is a del Pezzo surface of degree $2,3,4,5,6,7,6$ respectively, $H_{1i}$ are $K3$ surfaces anticanonical in $Y_i$, $C = D \bigcap H_{1i}$ is a smooth curve of genus $3,4,5,6,7,8,7$ respectively and the ample polarization $H_{11}' \bigcup_{C'} H_{12}'$ consists of $K3$ surfaces $H_{1i}'$ anticanonical in $Y_i$ with $H_{1i} \neq H_{1j}'$ and a smooth curve $C' = D \bigcap H_{1i}'$ of genus $3,4,5,6,7,8,7$ respectively with $C \neq C'$.
        
        \item[(b)] In the general fibre, we have a pair $(V_t, \displaystyle\frac{1}{2}B_t+\epsilon B_t')$ where $V_t$ is Fano of type $1.2, 1.3, 1.4$, $1.5(b)$, $2.6(b)$, $2.8$ or $3.1$ respectively as mentioned in Corollary \ref{Fano-threefolds} while $B_t$ and the polarization $ B_t'$ both $ \in |-2K_{V_t}|$ are smooth surfaces with $B_t \neq B_t'$. 
    \end{itemize}

    \item[(2)] For $m = 2, n = 1, d = 4$, in cases $(1)-(4), (6)-(8)$ 
    \begin{itemize}
       
        \item[(a)]  in the central fibre we have a polarized Calabi-Yau pair $(Y_1 \bigcup_D Y_2, \displaystyle\frac{1}{4}(H_{11} \bigcup_C H_{12}) + \epsilon (H_{11}' \bigcup_{C'} H_{12}'))$ for $\epsilon < 1$, where $Y_i$, $D$, $H_{1i}'$ and $C'$  are same as in part $(1)(a)$, $H_{1i} \in |-2K_{Y_i}| $ are canonically embedded surfaces in $Y_i$, $C = D \bigcap H_{1i}$ is a curve of genus $13, 19, 25, 31, 37, 43, 37$ respectively 
        \item[(b)] In the general fibre, we have a pair $(V_t, \displaystyle\frac{1}{4}B_t+\epsilon B_t') \subset \mathbb{P}^N$ where $V_t$ and $B_t'$ are as in part $(1)(b)$ while $B_t \in |-4K_{V_t}|$ is a smooth surfaces
    \end{itemize}
    
     \item[(3)] For $m = 1, n = 2, d = 3$, in case $(5)(a)$
     \begin{itemize}
        
        \item[(a)] in the central fibre we have a polarized Calabi-Yau pair $(Y_1 \bigcup_D Y_2, \displaystyle\frac{2}{3}(H_{11} \bigcup_C H_{12})+ \epsilon (H_{11}' \bigcup_{C'} H_{12}'))$ for $\epsilon < 1$, where $Y_i$ are Fano threefolds of type $1.16$, $D \in -\displaystyle\frac{1}{3}K_Y$ is a del-Pezzo surface of degree $8$, $H_{1i}$ are $K3$ surfaces anticanonical in $Y_i$ and $C = D \bigcap H_{1i}$ is a smooth curve of genus $4$ and the ample polarization $H_{11}' \bigcup_{C'} H_{12}'$ consists of $K3$ surfaces $H_{1i}'$ anticanonical in $Y_i$ with $H_{1i} \neq H_{1j}'$ and $C' = D \bigcap H_{1i}'$ is a smooth curve of genus $4$ with $C \neq C'$.
        \item[(b)] In the general fibre, we have a pair $(V_t, \displaystyle\frac{2}{3}B_t+\epsilon B_t')$ where $V_t$ is Fano of type $1.5(b)$ as mentioned in Corollary \ref{Fano-threefolds}  while $B_t$ and the polarization $B_t'$ both $\in |-\displaystyle\frac{3}{2}K_{V_t}|$ are smooth surfaces with $B_t \neq B_t'$.
    \end{itemize}

    \item[(4)] For $m = 1, n = 1, d = 3$, in case $(5)(b)$
     \begin{itemize}
        
        \item[(a)] in the central fibre we have a polarized Calabi-Yau pair $(Y_1 \bigcup_D Y_2, \displaystyle\frac{1}{3}(H_{11} \bigcup_C H_{12})+ \epsilon (H_{11}' \bigcup_{C'} H_{12}'))$ for $\epsilon < 1$, where $Y_i$ are Fano threefolds of type $1.16$, $D \in -\displaystyle\frac{2}{3}K_Y$ is a del-Pezzo surface of degree $4$, $H_{1i}$ are $K3$ surfaces anticanonical in $Y_i$ and $C = D \bigcap H_{1i}$ is a smooth curve of genus $13$ and the ample polarization $H_{11}' \bigcup_{C'} H_{12}'$ consists of $K3$ surfaces $H_{1i}'$ anticanonical in $Y_i$ with $H_{1i} \neq H_{1j}'$ and $C' = D \bigcap H_{1i}'$ is a smooth curve of genus $13$ with $C \neq C'$.
        \item[(b)] In the general fibre, we have a pair $(V_t, \displaystyle\frac{1}{3}B_t+\epsilon B_t')$ where $V_t$ is Fano of type $1.2(b)$ as mentioned in Corollary \ref{Fano-threefolds}  while $B_t$ and the polarization $B_t'$ both $\in |-3K_{V_t}|$ are smooth surfaces with $B_t \neq B_t'$.
    \end{itemize}

  \end{itemize}
    
\end{example}

In the following Corollary we give a list of stable pairs $(V,\Delta)$ where $V = Y_1 \bigcup_D Y_2$ is a union of two Fano threefolds intersecting along a $K3$ surface $D$ and the pair is smoothable into stable pairs $(V_t, \Delta_t)$ where $V_t$ is smooth. 
\begin{corollary}\label{log gen type pairs from CY and gen tye}
The following smoothable stable pairs are obtained from Corollary \ref{Calabi-Yau threefolds} and Corollary \ref{Fano with Bott} using Theorem \ref{pairs main} :
Let $d, m$ and $n$, be integers that satisfy $d \geq n$. Then for each smoothable $V \subset \mathbb{P}^N$, in Corollary \ref{Calabi-Yau threefolds} (resp. Corollary \ref{Fano with Bott}), and a general hyperplane section $B_m \in \mathcal{O}_V(m)$ and there exists a Gorenstein one parameter family $(\mathcal{V}, \displaystyle\frac{n}{d}\mathcal{B}_m) \to T$ of pairs such that
    \begin{itemize}
        \item[(a)] Both $K_{\mathcal{V}/T}$ and $\mathcal{B}_m$ are Cartier.
        \item[(b)] $(\mathcal{V}_0, \displaystyle\frac{n}{d}\mathcal{B}_{m0}) = (V, \displaystyle\frac{n}{d}B_m)$ is a stable pair, where $B_m = H_{m1} \bigcup_C H_{m2}$ is a simple normal crossing Cartier divisor in $V$, such that $H_{mi} \in |\mathcal{O}_{Y_i}(m)|$ are smooth divisors in $Y_i$ and $C = H_{m1} \bigcap D = H_{m2} \bigcap D$ is smooth.
        \item[(c)] For $t \neq 0$,  $(\mathcal{V}_t, \displaystyle\frac{n}{d}\mathcal{B}_t) = (V_t, \displaystyle\frac{n}{d}B_{mt})$ is a stable lc pair with, where $V_t$ is a smooth Calabi-Yau threefold (resp. a threefold of general type) and $B_{mt} \subset V_t$ is a smooth Cartier divisor.
   \end{itemize}

\end{corollary}

\begin{example}\label{examples of log pairs from CY and gen type}
We describe explicitly some of the smoothable stable pairs obtained from Corollary \ref{log CY and gen type pairs from Fano}. 
\begin{itemize}
    \item[(1)] For $m = 1, n = 1, d = 1$ in Corollary \ref{Calabi-Yau threefolds},  
    \begin{itemize}
        
        \item[(a)]  in the central fibre we have a a stable pair $(Y_1 \bigcup_D Y_2, H_{11} \bigcup_C H_{12})$ where $Y_i$ are Fano threefolds both of same deformation type from one of families $1.5-1.10$ and $1.12-1.17$, $D$ is a polarized $K3$ surface anticanonical in $Y_i$, $H_{1i} \neq D$ are polarized $K3$ surfaces anticanonical in $Y_i$ and $C = D \bigcap H_{1i}$ is a canonically embedded smooth curve of appropriate genus.
        \item[(b)] In the general fibre, we have a pair $(V_t, B_t)$ where $V_t$ is a smooth Calabi-Yau threefold and $B_t$ is a canonically embedded surface. 
        
    \end{itemize}
   
    \item[(2)] For $m = 1, n = 1, d = 1$ in Corollary \ref{Fano with Bott},
    \begin{itemize}
        
        \item[(a)]  in the central fibre we have a stable pair $(Y_1 \bigcup_D Y_2, H_{11} \bigcup_C H_{12})$, where $Y_i$'s are Fano threefolds satisfying Bott vanishing theorem, both of same deformation type, $D \in |-2K_{Y_i}|$ is a smooth surface, $H_{1i} \in |-2K_{Y_i}|$ are smooth surfaces in $Y_i$ and $C = D \bigcap H_{1i}$ is a smooth curve of appropriate genus.
        
        \item[(b)] In the general fibre, we have a pair $(V_t, B_t)$ where $V_t$ is a threefold with ample canonical bundle and $B_t \in |2K_{V_t}|$ is a smooth surface.
    \end{itemize}
     
\end{itemize}
    
\end{example}

\vspace{0.5cm}

\bibliographystyle{plain}

\end{document}